\documentclass[10pt,reqno]{amsart}
\usepackage[all]{xy}
\usepackage{amscd}
\usepackage{xcolor}
\usepackage{a4wide}
\usepackage{amssymb}
\usepackage{amsfonts}
\usepackage{amsmath}
\input xy
\xyoption{arrow} \xyoption{matrix}

\date{}

\newtheorem{proposition}{Proposition}[section]
\newtheorem{theorem}[proposition]{Theorem}
\newtheorem{lemma}[proposition]{Lemma}
\newtheorem{example}[proposition]{Example}

\newtheorem{corollary}[proposition]{Corollary}

\def\Supp{\mathop\mathrm{Supp}}

\def\der{\partial }

\def\nFM0{{\nu }_{F,M_0}}
\def\nFN0{{\nu }_{F,N_0}}
\def\nGN0{{\nu }_{G,N_0}}

\def\N0{ {\bf N}_0 }

\def\t{\otimes}

\def\ra{\rightarrow}

\def\Xpm{X^{\pm }}

\def\s{\sigma}
\def\Z{\mathbb{Z}}

\def\l1{{\lambda}_1}

\def\a{\alpha}
\def\a0{ {\alpha }_0}
\def\a1{ {\alpha }_1}

\def\cG{{\mathcal G}}
\def\cM{{\mathcal M}}
\def\l{\lambda}

\def\nFGM0{{\nu }_{F,G,M_0}}

\def\nFN0{{\nu}_{F,N_0}}

\def\sm{{\sigma}^m}

\def\sm1{{\sigma}^{-1}}

\def\smtp1{{\sigma}^{-t+1}}

\def\S1{S^{-1}}

\def\Xpm1{X^{\pm 1}_1}

\def\sPM1{{\sigma }^{\pm 1}}
\def\sMP1{{\sigma }^{\mp 1 }}

\def\d{\delta}

\def\di{{\rm d.ind}}

\def\L{\Lambda}
\def\O{\Omega}

\def\G{\Gamma}

\def\OO{{\mathcal O}}

\def\CD{{\mathcal D}}

\def\Ytm1{Y^{t-1}}
\def\Yim1{Y^{i-1}}

\def\CN{{\mathcal N}}

\def\Aut{{\rm Aut}}

\def\dim{{\rm dim }}

\def\SL2Z{ {\rm SL}_2({\bf Z}) }

\def\th{ \theta }

\def\Gp1{ G^{1 , 1 } }
\def\P11{ P^{-1 , 1 } }
\def\Pp1{ P^{1 , 1 } }

\def\Supp{{\rm Supp}}

\def\th{\theta}

\def\nCLsr{{}^\nu\kern-2pt {\cal L}^{\sigma , \rho  }}
\def\nP{{}^\nu \kern-2pt P}
\def\nL{{}^\nu\kern-2pt L}
\def\nLL{{}^\nu\kern-2pt \Lambda}
\def\nPsr{{}^\nu\kern-2pt P^{\sigma , \rho  }}
\def\nLsr{{}^\nu\kern-2pt L^{\sigma , \rho  }}
\def\nuCL{{}^\nu\kern-2pt  {\cal L}}
\def\nCLsr{{}^\nu\kern-2pt {\cal L}^{\sigma , \rho  }}
\def\nCL1m{{}^\nu\kern-2pt {\cal L}^{-1 , 1  }}
\def\x1nu{x^\frac{1}{\nu}}
\def\xm1nu{x^{-\frac{1}{\nu}}}

\def\CN{{\mathcal N}}
\def\ra{\rightarrow }

\def\CB{{\cal B}}

\def\CI{{\mathcal I}}

\def\CC{ {\mathcal C}}
\def\O{{\mathcal O}}

\def\nAM0{{\nu }_{{\cal A},M_0}}
\def\nAN0{{\nu }_{{\cal A},N_0}}

\def\End{ {\rm End }}
\def\Max{ {\rm Max }}

\def\ga{\mathfrak{a}}
\def\gb{\mathfrak{b}}

\def\gn{\mathfrak{n}}
\def\gm{\mathfrak{m}}
\def\gp{\mathfrak{p}}

\def\SL{{\rm SL}}

\def\Spec{{\rm Spec}}

\def\Ext{{\rm Ext}}

\def\di!{\frac{\der^i}{i!}}
\def\dik!{\frac{\der^k_i}{k!}}

\def\id{{\rm id}}

\def\Md#1{#1\mbox{-}\mathrm{Mod}}

\def\Max{{\rm Max}}

\def\N{\mathbb{N}}

\def\0{\overline{0}}
\def\1{\overline{1}}

\def\Ln1{\L_{n,\overline{1}}}

\def\a1{a_{\overline{1}}}

\def\S{\Sigma}

\def\vn1{\overrightarrow{n-1}}

\def\sl{{\rm sl}}

\def\mA{\mathbb{A}}
\def\mD{\mathbb{D}}

\def\soc{{\rm soc}}

\def\mJ{\mathbb{J}}
\def\mI{\mathbb{I}}

\def\ann{{\rm ann}}

\def\bM{\overline{M}}

\def\ind{{\rm ind}}

\def\K1{{\rm K}_1}

\def\hmI1{\widehat{\mI_1}}
\def\tmI1{\widetilde{\mI_1}}
\def\tmJ1{\widetilde{\mJ_1}}
\def\hB1{\widehat{B_1}}
\def\hCB1{\widehat{\CB_1}}

\def\WIn{{\rm W(\mI_n)}}
\def\WO{{\rm W(\mI_n, \O)}}
\def\WDO{{\rm W(\mI_n, \CD_\O)}}
\def\GW{{\rm GW}}
\def\GWIn{{\rm GW(\mI_n)}}
\def\GWO{{\rm GW(\mI_n, \O)}}
\def\GWDO{{\rm GW(\mI_n, \CD_\O)}}
\def\rGWDO{{\rm GW_r(\mI_n, \CD_\O)}}
\def\rWDO{{\rm W_r(\mI_n, \CD_\O)}}
\def\GWDOi{{\rm GW(\mI_n, \CD_{\O_{i}})}}
\def\Ker{\mathop\mathrm{Ker}\nolimits}
\def\LF{{\rm LF}}
\def\LFm{{\LF_{\gm}}}
\def\Od{{\O_{\rm deg}}}
\def\Ond{{\O_{\rm ndeg}}}

\def\rWIn{{\rm W_r(\mI_n)}}
\def\rWO{{\rm W_r(\mI_n, \O)}}
\def\rGW{{\rm GW_r}}
\def\rGWIn{{\rm GW_r(\mI_n)}}
\def\rGWO{{\rm GW_r(\mI_n, \O)}}
\def\OOO{\mathbb{O}}

\begin{document}

\author{V. V.  \  Bavula, V.  Bekkert  \and V. Futorny}
\title[Indecomposable weight  modules]{EXPLICIT DESCRIPTION OF GENERALIZED WEIGHT MODULES OF THE ALGEBRA OF POLYNOMIAL INTEGRO-DIFFERENTIAL OPERATORS $\mI_n$}
\address{Department of Pure Mathematics, University of Sheffield, Hicks Building, Sheffield S3 7RH, UK}
\email{v.bavula@sheffield.ac.uk}
\address{Departamento de Matem\'atica, ICEx, Universidade Federal de Minas
Gerais, Av.  Ant\^onio Carlos, 6627, CP 702, CEP 30123-970, Belo
Horizonte-MG, Brasil} \email{bekkert@mat.ufmg.br}
\address{Instituto de Matem\'atica e Estat\'{\i}stica,
Universidade de S\~ao Paulo, Caixa Postal 66281, S\~ao Paulo, CEP
05315-970, Brasil} \email{futorny@ime.usp.br}



\begin{abstract}

For the algebra $\mI_n= K\langle x_1, \ldots, x_n, \der_1, \ldots, \der_n, \int_1, \ldots, \int_n \rangle$ of
polynomial integro-differential operators over a field $K$ of
characteristic zero, a classification of simple weight and generalized weight (left and right) $\mI_n$-modules   is given. 
It is proven that the category of  weight $\mI_n$-modules is semisimple. An explicit 
description of generalized weight $\mI_n$-modules is given and using it a criterion is obtained for the problem
of classification of indecomposable generalized weight $\mI_n$-modules to be of finite representation type, tame or wild.
In the tame case, a classification of indecomposable generalized weight $\mI_n$-modules   is given. 
In the wild case `natural` tame subcategories are considered with explicit description of  indecomposable modules.
It is proven that every generalized weight $\mI_n$-module
is a unique sum of absolutely prime modules.
For an arbitrary ring $R$, we introduce the concept of
{\em absolutely prime} $R$-module (a nonzero $R$-module $M$
is absolutely prime if all nonzero subfactors of $M$ have the
same annihilator).
It is shown that every indecomposable generalized weight $\mI_n$-module is equidimensional. 
A criterion is given for a generalized weight $\mI_n$-module to be finitely generated.




\end{abstract}

\maketitle


\section{Introduction}

Throughout, ring means an associative ring with $1$; module means
a left module;
 $\N :=\{0, 1, \ldots \}$ is the set of natural numbers; $\N_+ :=\{ 1,2,  \ldots \}$ and $\Z_{\leq 0} :=-\N $;   $K$ is a
field of characteristic zero and  $K^*$ is its group of units;
$\otimes =\otimes_K$; $P_n:= K[x_1, \ldots , x_n]$ is a polynomial algebra over $K$;
$\der_1:=\frac{\der}{\der x_1}, \ldots , \der_n:=\frac{\der}{\der
x_n}$ are the partial derivatives ($K$-linear derivations) of
$P_n$; $\End_K(P_n)$ is the algebra of all $K$-linear maps from
$P_n$ to $P_n$;
the subalgebras  
$$A_n:= K \langle x_1, \ldots , x_n , \der_1,
\ldots , \der_n\rangle \;\;{\rm and} \;\;\mI_n:=K\langle x_1, \ldots , x_n,
 \der_1, \ldots ,\der_n,  \int_1,
\ldots , \int_n\rangle $$  
\noindent of the algebra $\End_K(P_n)$ are called
the $n$'th {\em Weyl algebra} and the {\em algebra of  polynomial
integro-differential operators}, respectively.

The Weyl algebras $A_n$ are Noetherian algebras and domains. The
algebras $\mI_n$  were introduced in \cite{algintdif,intdifaut} 
 and they  are neither left nor right Noetherian and not
domains. Moreover, they contain infinite direct sums of nonzero
left and right ideals \cite{algintdif}. 
The algebra $\mI_n$ contains a polynomial algebra 
$$D_n=K[H_1,\ldots , H_n], \;\; {\rm where}\;\;
H_1:=\der_1x_1,\ldots , H_n:=\der_nx_n,$$ 
\noindent which is a maximal commutative subalgebra of $\mI_n$,
\cite{algintdif}. 
An $\mI_n$-module $M$ is called a {\em weight module} if it is a {\em semisimple} $D_n$-module provided the field $K$
is an algebraically closed field.
For an arbitrary field, a weight $\mI_n$-module is a direct sum of  common eigenspaces for the commuting elements $H_1,\ldots,H_n$, by definition.
 An $\mI_n$-module $M$ is called a {\em generalized weight module} if for all  elements $m\in M$, $\dim_K(D_nm)<\infty$, i.e., $M$
 is a {\em locally finite dimensional} $D$-module provided the field $K$
is an algebraically closed field. 
 Every weight module is a generalized weight module but not vice versa.
 The $\mI_n$-module $P_n$ is a simple weight module.
 Introduction of  (generalized) 
 weight modules for the algebras $\mI_n$ was inspired by a similar concept  for semisimple finite dimensional Lie algebras.
 For a semisimple finite dimensional  Lie algebra  $\cG$,  a  classification of {\em simple} generalized  weight modules is known only when $\cG=\sl_2$.
 Furthermore, a classification of {\em indecomposable}  generalized weight $\sl_2$-modules was done in \cite{Bav-Bek} and 
 the problem of classification turned out to be tame provided the Casimir element acts as a scalar. In fact, in \cite{Bav-Bek}  in the tame case
 a classification of indecomposable  
 generalized weight modules was obtained for a large class of algebras, the, so-called, {\em generalized Weyl algebras} 
 (the universal enveloping algebra $U(\sl_2)$ and the Weyl algebras are examples of generalized Weyl algebras as well as many other quantum groups are). Recently, for some non-semisimple Lie algebras $\cG$ and their quantum analogues, classifications of simple weight modules are given: the Schr\"{o}dinger algebra, \cite{Bav-Lu-N9}; $\sl_2\ltimes V_2$,\;\cite{Bav-Lu-N5}, where $V_2$ is the simple 2-dimensional  $\sl_2$-module; the enveloping algebra of the Euclidean algebra,  \cite{Bav-Lu-N8};  $\mathbb{K}_q[X,Y]\rtimes U_q(\sl_2)$, \cite{Bav-Lu-N7}; the quantum spatial ageing algebra, \cite{Bav-Lu-N3}.

 In the paper, explicit descriptions of weight and generalized  weight $\mI_n$-modules are given (Theorem \ref{15May19}, Theorem \ref{14Apr17} and (\ref{MMO2}), (\ref{GWDD})). 
 They are too technical to explain them in the Introduction.   Classifications of simple weight and simple generalized weight $\mI_n$-modules are obtained (Theorem \ref{A11Apr17}). 
 It is proven that the category of weight $\mI_n$-modules is a semisimple category (Theorem \ref{B11Apr17}).
 In \cite{BBF-1}, a classification of indecomposable generalized weight $\mI_1$-modules of finite length is given. 
This classification is used in the present paper in order to obtain the general case.
Usually, the case $n=1$ serves as the base of induction for the case $n>1$.
 Using Theorem \ref{A11Apr17}, a criterion is given to decide whether the problem of classification of indecomposable generalized 
weight $\mI_n$-modules is of finite representation type, tame or wild (Theorem \ref{A14Apr17}).
In the case, $n=1$ the problem is tame, \cite{BBF-1}. 
In the tame case, a classification of indecomposable generalized weight $\mI_n$-modules   is given.
It is shown that every indecomposable generalized weight $\mI_n$-module is equidimensional (Corollary \ref{b15Apr17}). 
A criterion is given for a generalized weight $\mI_n$-module to be finitely generated  (Corollary \ref{c15Apr17}).
In the wild case,   `natural' tame subcategories are considered with explicit description of  indecomposable modules,
see the end of Section \ref{CSGWM}. In particular, descriptions of categories $\ind(D_2, \gm^{2})$, $\ind_f(\G)$ and $\ind_f(A)$ are
obtained.
In Section \ref{GWRIM}, similar results are proven for generalized weight {\em right} $\mI_n$-modules.

 Properties of the algebras $\mI_n$  are studied in \cite{algintdif, algintdifline, indtif-bimod}. In the case $n=1$, for a more general setting   see also \cite{Guo-Regen-Ros-2014}. 
 The simple $\mI_1$-modules are classified in \cite{algintdifline}. Simple $A_1$-modules were classified in \cite{Block-IrrRepsl2} (see also \cite{Bav-GWA-1992, Bav-92}  
 for some generalized Weyl algebras including $A_1$). The automorphism groups $\Aut_{K-{\rm alg}}(\mI_n)$ are found in \cite{intdifaut}.
 The weak homological dimension of the algebra $\mI_n$ is $n$, \cite{algintdif}. Futhermore, the weak homological dimension is $n$ for all the factor algebras of $\mI_n$, \cite{Bav-gldim-intdif-JacAlg}.

Finite dimensionality of Ext-groups of simple modules over the (first) Weyl algebra $A_1$  was proven in \cite{MR}. 
Finite dimensionality of Ext-groups of simple modules over the generalized  Weyl algebras   was proven in \cite{Bav-FinDimExt-91}. 
Simple modules over certain generalized Weyl algebras were classified in \cite{Bav-GWA-1992}.
In \cite{BBF-1}, Ext-groups are described between indecomposable generalized weight $\mI_1$-modules, it is shown  that they are finite dimensional vector spaces.
In \cite{indtif-bimod}, it is proven that the algebra $\mI_n$ is a left coherent algebra iff the algebra $\mI_n$ is a right coherent iff $n=1$;
the algebra $\mI_n$ is a maximal left (resp., right) order in the largest left (resp., right) quotient ring $Q_l(\mI_n)$ (resp., $Q_r(\mI_n)$) of $\mI_n$.
The (left and right)  global dimension   of the algebra $\mI_n$  and all prime factor  algebras of $\mI_n$ is equal to $n$, \cite{Bav-gldim-intdif-JacAlg}.

Classifications of (various classes of) simple weight modules over algebras that are close to the (generalized) Weyl algebras  
are given in \cite{Maz-Tur-1999,Bek-Ben-Fut-2004,Hartwig-2006,Ship-2010,Hartwig-2011,Fut-Gr-Maz-2014,Lu-Maz-Zha0-2015,Fut-Iyer-2016,Bav-Lu-2016-Isr}.


\section{Classification of simple (generalized) weight  $\mI_n$-modules}\label{CSGWM}


In this section, a classification  of simple generalized weight and simple weight  $\mI_n$-modules  is
given (Theorem \ref{A11Apr17}). It is proven that the category of  weight $\mI_n$-modules
is a semisimple category (Theorem \ref{B11Apr17}). At the beginning of the section, we collect some
results about the algebras $\mI_n$ that are used in the paper. In the case when $n=1$, we drop the subscript `$1$'  
in order to simplify the notation.

As an abstract algebra, the algebra $\mI_1$  is generated by the elements $\der $, $H:=
\der x$ and $\int$ (since $x=\int H$) that satisfy the defining
relations, {\cite[Proposition 2.2]{algintdif}} (where $[a,b]:=ab-ba$): 

\begin{equation}\label{DefRI}
\der \int = 1,
\;\; [H, \int ] = \int, \;\; [H, \der ] =-\der , \;\; H(1-\int\der
) =(1-\int\der ) H = 1-\int\der .
\end{equation}

Since $\mI_n=\mI_1\otimes \cdots \otimes \mI_1$ ($n$ times), defining relations  of the algebra $\mI_n$ is the union of the 
defining relations (\ref{DefRI}) for each index $i=1, \ldots, n$ and the relations $a_ia_j=a_ja_i$
for all $i\neq j$ where $a_i\in \{\der_i, H_i, \int_i\}$.
The elements of the algebra
$\mI_1$, 
\begin{equation}\label{eijdef}
e_{ij}:=\int^i\der^j-\int^{i+1}\der^{j+1}, \;\; i,j\in \N ,
\end{equation}
satisfy the relations $e_{ij}e_{kl}=\d_{jk}e_{il}$ where $\d_{jk}$
is the Kronecker delta function. Notice that
$e_{ij}=\int^ie_{00}\der^j$. The matrices of the linear maps
$e_{ij}\in \End_K(K[x])$ with respect to the basis $\{ x^{[s]}:=
\frac{x^s}{s!}\}_{s\in \N}$ of the polynomial algebra $K[x]$  are
the elementary matrices, i.e.
$$ e_{ij}*x^{[s]}=\begin{cases}
x^{[i]}& \text{if }j=s,\\
0& \text{if }j\neq s.\\
\end{cases}$$
Let $E_{ij}\in \End_K(K[x])$ be the usual matrix units, i.e.
$E_{ij}*x^s= \d_{js}x^i$ for all $i,j,s\in \N$. Then
\begin{equation}\label{eijEij}
e_{ij}=\frac{j!}{i!}E_{ij},
\end{equation}
 $Ke_{ij}=KE_{ij}$, and
$F:=\bigoplus_{i,j\geq 0}Ke_{ij}= \bigoplus_{i,j\geq
0}KE_{ij}\simeq M_\infty (K)$, the algebra (without 1) of infinite
dimensional matrices.
The algebra $\mI_n=\mI_1(1)\t\cdots \t\mI_1(n)\simeq \mI_1^{\t n}$
where $\mI_1(i)=K\langle  \der_i, H_i, \int_i\rangle $ for $i=1, \ldots, n$.\\

{\bf $\Z^{n}$-grading on the algebra $\mI_n$ and the canonical form of
an integro-differential operator, \cite{algintdif}.} The algebra
$\mI_1=\bigoplus_{i\in \Z} \mI_{1, i}$ is a $\Z$-graded algebra
($\mI_{1, i} \mI_{1, j}\subseteq \mI_{1, i+j}$ for all $i,j\in
\Z$) where  
\begin{equation}\label{I1iZ}
\mI_{1, i} =\begin{cases}
D_1\int^i=\int^iD_1& \text{if } i>0,\\
D_1& \text{if }i=0,\\
\der^{|i|}D_1=D_1\der^{|i|}& \text{if }i<0,\\
\end{cases}
\end{equation}
 the algebra $D_1:= K[H]\bigoplus \bigoplus_{i\in \N} Ke_{ii}$ is
a commutative non-Noetherian subalgebra of $\mI_1$, $ He_{ii} =
e_{ii}H= (i+1)e_{ii}$  for $i\in \N $ (notice that
$\bigoplus_{i\in \N} Ke_{ii}$ is the direct  sum of non-zero
ideals of $D_1$); 
$$(\int^iD_1)_{D_1}\simeq D_1,\;\; \int^id\mapsto
d;\;\; {}_{D_1}(D_1\der^i) \simeq D_1,\;\; d\der^i\mapsto d \;\;{\rm  for
\;all}\;\; i\geq 0$$ 
\noindent since $\der^i\int^i=1$.
 Notice that the maps $\cdot\int^i : D_1\ra D_1\int^i$, $d\mapsto
d\int^i$,  and $\der^i \cdot : D_1\ra \der^iD_1$, $d\mapsto
\der^id$, have the same kernel $\bigoplus_{j=0}^{i-1}Ke_{jj}$.
The algebra 
$$\mI_n=\bigoplus_{\alpha\in\Z^{n}}\mI_{n,\alpha}$$ 
\noindent is a $\Z^{n}$-graded algebra
($\mI_{n,\alpha}\mI_{n,\beta}\subseteq \mI_{n,\alpha+\beta}$ for all $\alpha, \beta\in \Z^{n}$) 
where $\mI_{n,\alpha}=\otimes_{i=1}^{n}\mI_1(i)_{\alpha_i}$ and $\alpha=(\alpha_1, \ldots, \alpha_n)$.

Each element $a$ of the algebra $\mI_1$ is a unique finite sum
\begin{equation}\label{acan}
a=\sum_{i>0} a_{-i}\der^i+a_0+\sum_{i>0}\int^ia_i +\sum_{i,j\in
\N} \l_{ij} e_{ij}
\end{equation}
where $a_k\in K[H]$ and $\l_{ij}\in K$. This is the {\em canonical
form} of the polynomial integro-differential operator, 
\cite{algintdif}. Let $$v_i:=\begin{cases}
\int^i& \text{if }i>0,\\
1& \text{if }i=0,\\
\der^{|i|}& \text{if }i<0.\\
\end{cases}$$

\noindent Then $\mI_{1,i}=D_1v_i= v_iD_1$ and an element $a\in \mI_1$ is the
unique  finite  sum 
\begin{equation}\label{acan1}
a=\sum_{i\in \Z} b_iv_i +\sum_{i,j\in \N} \l_{ij} e_{ij}
\end{equation}
where $b_i\in K[H]$ and $\l_{ij}\in K$. So, the set $\{ H^j\der^i,
H^j, \int^iH^j, e_{st}\, | \, i\geq 1; j,s,t\geq 0\}$ is a
$K$-basis for the algebra $\mI_1$. The tensor product of these bases is a basis for the algebra $\mI_n$.
The multiplication in the
algebra $\mI_1$ is given by the rule:
$$ \int H = (H-1) \int , \;\; H\der = \der (H-1), \;\; \int e_{ij}
= e_{i+1, j},$$
$$ e_{ij}\int= e_{i,j-1}, \;\; \der e_{ij}=
e_{i-1, j},\;\; e_{ij} \der = \der e_{i, j+1},$$
$$ He_{ii} = e_{ii}H= (i+1)e_{ii}, \;\; i\in \N, $$
where $e_{-1,j}:=0$ and $e_{i,-1}:=0$.  $\noindent $

 The algebra
$\mI_1$ has the only proper ideal $F=\bigoplus_{i,j\in \N}Ke_{ij}
\simeq M_\infty (K)$ and
 $F^2= F$. The factor algebra $\mI_1/F$ is canonically isomorphic to
the skew Laurent polynomial algebra $$B_1:= K[H][\der, \der^{-1} ;
\tau ],\;\; \tau (H) = H+1, \;\; {\rm via}\;\; \der \mapsto \der,\;\;  \int\mapsto
\der^{-1},\;\; H\mapsto H$$ 
\noindent (where $\der^{\pm 1}\alpha = \tau^{\pm
1}(\alpha ) \der^{\pm 1}$ for all elements $\alpha \in K[H]$). The
algebra $B_1$ is canonically isomorphic to the (left and right)
localization $A_{1, \der }$ of the Weyl algebra $A_1$ at the
powers of the element $\der$ (notice that $x= \der^{-1} H$).

Recall that the algebra of polynomial integro-differential operators 
$$\mI_n=K\langle x_1, \ldots, x_n, \der_1, \ldots, \der_n, \int_1, \ldots, \int_n\rangle$$
\noindent over
a field $K$ of characteristic zero $\mI_n=\mI_1(1)\otimes \cdots\otimes \mI_1(n)$ is the tensor product
of algebras $\mI_1(i)=K\langle x_i, \der_i, \int_i\rangle\simeq \mI_1$.
Each algebra $\mI_1(i)$ contains a unique proper ideal 
$$F(i)=\bigoplus_{s,t\in \N}Ke_{st}(i) \;\;
{\rm where}\;\; e_{st}=\int_i^{s}\der_i^{t}-\int_i^{s+1}\der_i^{t+1},$$ 
\noindent $F(i)^{2}=F(i)$
and $B_1(i):=\mI(i)/F(i)\simeq K[H_i][\der_i, \der_i^{-1}; \s_i^{-1}]$ where $\s_i(H_i)=H_i-1$.
The algebra $\mI_n$ is a local algebra where the unique maximal ideal $\ga_n$ is generated by
$F(1), \ldots, F(n)$ and the factor algebra $\mI_n/\ga_n$ is isomorphic to the skew Laurent polynomial
algebra 
$$B_n=D_n[\der_1^{\pm 1}, \ldots, \der_n^{\pm 1}; \s_1^{-1}, \ldots, \s_n^{-1}]\;\;{\rm  where}\;\;
\s_i(H_j)=H_j -\delta_{ij},$$ 
\noindent see \cite{algintdif}. 
Furthermore, the algebra $B_n$ is the only left/right Noetherian factor algebra of $\mI_n$, \cite{algintdif}.

 A classification of all ideals (including prime
ideals) of the algebra $\mI_n$ is obtained in \cite{algintdif}.
There are precisely $n$ height $1$ prime ideals:
$$ \gp_1=F\t \mI_{n-1}, \;\; \gp_2=\mI_1\t F\t \mI_{n-2}, \ldots, \gp_n = \mI_{n-1}\t F,$$
see \cite{algintdif}. The algebra $\mI_n$ is a prime algebra ($0$ is a prime ideal of $\mI_n$). Every nonzero prime ideal $\gp$ is a unique sum 
$$\gp_I=\sum_{i\in I}\gp_i$$ 
\noindent of height 1 prime ideals where $I\subseteq \{ 1, \ldots , n\}$ is a unique set for the ideal $\gp = \gp_I$, and ${\rm ht} (\gp )=|I|$ where ${\rm ht}(\gp ) $ if the height of the ideal $\gp$. Every ideal of $\mI_n$ is an idempotent ideal ($\ga^2=\ga$), ideals of $\mI_n$ commute ($\ga \gb =\gb\ga$) and the ideal $\ga_n=\gp_1+\cdots +\gp_n$ is the only maximal ideal of the algebra $\mI_n$.\\

{\bf Generalized weight $\mI_n$-modules.}
The group $\Z^{n}$ acts on the vector space $K^{n}$ by addition.
For an element $\l=(\l_1, \ldots , \l_n)\in K^{n}$, $\O(\l)=\l +\Z^{n}$ is its orbit.
The set of all $\Z^{n}$-orbits is isomorphic to the factor group $K^{n}/\Z^{n}$,
$\O(\l)\leftrightarrow \lambda +\Z^{n}$.
In particular, two orbits are equal, $\O(\l)=\O(\l^{\prime})$  iff $\l-\l^{\prime}\in \Z^{n}$.

Let $e_1, \ldots, e_n$ be the standart basis for the vector space $K^{n}$. 
Then $\Z^{n}=\oplus_{i=1}^{n}\Z e_i\subseteq K^{n}=\oplus_{i=1}^{n}Ke_i$.

The $K$-automorphisms $\sigma_1,\ldots, \sigma_n$ of the polynomial algebra $D_n$ commute, $\sigma_i\sigma_j=\sigma_j\sigma_i$,
since $\sigma_i(H_j)=H_j-\delta_{ij}$ where $\delta_{ij}$ is the Kronecker delta.
The subgroup $G=\langle\sigma_1, \ldots, \sigma_n\rangle$ of the group of $K$-algebra automorphisms $\Aut_K(D_n)$ generated by the automorphisms
$\sigma_1,\ldots, \sigma_n$ is an abelian group isomorphic to $\Z^{n}$ via isomorphism $G\rightarrow \Z^{n}$, $\sigma_1\mapsto e_1$, $\ldots,$
$\sigma_n\mapsto e_n$. 

Let $\cM_n$ be the set of all maximal ideals of the algebra
$\mI_n$ of the type 
$$\gm =\gm_{\lambda}=(H_1-\lambda_1, \ldots,
H_n - \lambda_n)\;\; {\rm where}\;\; \lambda=(\lambda_1, \ldots, \lambda_n)\in K^{n}.$$ 
The group $\Aut_K(D_n)$ acts on $\cM_n$, $(\sigma,\gm) \mapsto \sigma(\gm)$.
Recall that the group $\Z^{n}$ acts on $K^{n}$ in the obvious way: $\Z^{n}\times K^{n}\rightarrow K^{n},$ $(i,\lambda)\mapsto i+\lambda$.
In a similar way, the group $G$ acts on $\cM_n$, $(\sigma_1^{i_1}\cdots\sigma_n^{i_n},\gm) \mapsto \sigma_1^{i_1}\cdots\sigma_n^{i_n}(\gm)$ where
$i=(i_1,\cdots, i_n)\in \Z^{n}$. 
When we identify the set $\cM_n$ with $K^{n}$ via the bijection 
$$\cM_n\rightarrow K^{n},\;\; \gm_{\lambda}\mapsto \lambda$$
\noindent and the group $G$ with $\Z^{n}$ via the group isomorphism $G\rightarrow \Z^{n}$, $\sigma_i\mapsto e_i$ $(i=1,\ldots, n)$, 
then the set of $G$-orbits $\cM_n/G$ is identified with the factor group $K^{n}/\Z^{n}$ via the bijection 
$$\cM_n/G\rightarrow K^{n}/\Z^{n},\;\;
G\gm_{\lambda}\mapsto \lambda+\Z^{n}$$
\noindent and the action $\sigma_1^{i_1}\cdots\sigma_n^{i_n}(\gm_{\lambda})$ is identified with the action $\lambda+i$ where $i=(i_1,\ldots, i_n)\in \Z^{n}$.

The polynomial algebra $D_n=K[H_1,\ldots,H_n]$ is a maximal commutative subalgebra of $\mI_n$.
For each  $\lambda=(\lambda_1,\ldots, \lambda_n)\in K^{n}$,
the field 
$$K_{\gm}:=D_n/\gm_{\lambda}=D_n/(H_1-\lambda_1,\ldots, H_n-\lambda_n)=K_{\lambda}\simeq K$$
\noindent is a unique simple
$D_n$-module that is annihilated by the maximal ideal $\gm_{\lambda}=(H_1-\lambda_1,\ldots, H_n-\lambda_n)$.
Let $\Max(D_n)$ be the maximum spectrum of $D_n$ (the set of maximal ideals of $D_n$).
Notice that if $K$ is an algebraically closed field then
$\cM_n=\Max(D_n)$.

An $\mI_n$-module $M$ is called a {\em weight} module if 
$$M=\bigoplus_{\gm\in \cM_n}M_{\gm}\;\; {\rm  where}\;\; M_{\gm}={\rm ann}_M(\gm):=\{m\in M\,|\, \gm m=0\}.$$
If $\gm=\gm_{\lambda}=(H_1-\lambda_1,\ldots, H_n-\lambda_n)$ for some $\lambda\in K^{n}$ then
$$M_{\lambda}:=M_{\gm_{\lambda}}=\{m\in M\,|\, H_1 m=\lambda_1 m, \ldots, H_n m=\lambda_n m\}.$$
The set $\Supp(M)=\{\gm\in \cM_n\, |\, M_{\gm}\neq 0\}$ is called the {\em support} of the weight
$\mI_n$-module $M$.
So, an $\mI_n$-module $M$ is weight iff it is a (direct) sum
of common eigen-spaces for the commuting elements $H_1, \ldots, H_n$
of the algebra $\mI_n$.

An $\mI_n$-module $M$ is called a {\em generalized weight} module if 
$$M=\bigoplus_{\gm\in \cM_n}M^{\gm}\;\; {\rm where}\;\;
M^{\gm}=\{m\in M\,|\, \gm^{i} m=0 \;\;{\rm for\; some}\;\; i\geq 0\}=\bigcup_{i\geq 0}{\rm ann}_M(\gm^{i}).$$
If $\gm=\gm_{\lambda}$ for some $\lambda\in K^{n}$ then 
$$M^{\gm}=\{m\in M\, |\, (H_1-\lambda_1)^{i}m=0,\ldots ,
(H_n-\lambda_n)^{i}m=0 \;\;{\rm for \; some}\;\; i\geq 1 \}.$$

The set $\Supp(M)=\{\gm\in \cM_n\, |\, M^{\gm}\neq 0\}$ is called the {\em support} of the generalized weight $\mI_n$-module $M$. 
Recall that we identified $(G, \cM_n, \cM_n/G)$ with
$(\Z^{n}, K^{n},  K^{n}/\Z^{n})$. 
Therefore, $\gm_{\lambda}$ is identified with $\lambda$,
$M_{\lambda}:=M_{\gm_{\lambda}}$ and $M^{\lambda}:=M^{\gm_{\lambda}}$. So, 
$$\Supp(M)=\{\lambda\in K^{n}\, |\, M^{\lambda}\neq 0\}.$$

If, in addition, the field $K$ is an algebraically closed then
the set $\Max(D_n)$ of maximal ideals of the ring  $D_n$ is 
equal to $\cM_n$, every weight $\mI_n$-module is an $\mI_n$-module
which is a semisimple $D_n$-module (and vice versa), and every
generalized weight $\mI_n$-module is an $\mI_n$-module which is 
{\em locally finite dimensional} $D_n$-module (and vice versa).

We denote by $\WIn$ (resp., $\GWIn$) the category of weight (resp., generalized weight) $\mI_n$-modules. Clearly,
$$\WIn\subseteq \GWIn\subseteq \Md{\mI_n}$$
are inclusions of categories where $\Md{\mI_n}$ is the category of all left $\mI_n$-modules.
The category $\GWIn$ is a full subcategory of $\Md{\mI_n}$, it is closed under arbitrary direct sums, extensions, submodules and factor modules.
The category $\WIn$ is closed under direct sums, submodules and factor modules but not under extensions, see  \cite[Theorem 2.5]{BBF-1}.

Let $M$ be a generalized weight $\mI_n$-module. It follows from the defining relations of the algebra $\mI_n$,  \cite[Proposition 2.2]{algintdif} or (\ref{DefRI}),   that for all
$\gm\in \cM_n$ and $i=1, \ldots, n$, 
$$x_iM^{\gm}\subseteq M^{\sigma_i(\gm)},\;\; \int_iM^{\gm}\subseteq M^{\sigma_i(\gm)},\;\; \der_iM^{\gm}\subseteq M^{\sigma_i^{-1}(\gm)}, $$
$$x_iM_{\gm}\subseteq M_{\sigma_i(\gm)},\;\; \int_iM_{\gm}\subseteq M_{\sigma_i(\gm)},\;\; \der_iM_{\gm}\subseteq M_{\sigma_i^{-1}(\gm)}. $$

So, the generalized weight $\mI_n$-module $M$ is a direct sum of its generalized weight submodules $M^{\O}$,

\begin{equation}\label{MMO}
M=\bigoplus_{\O\in\cM_n/G}M^{\O} \,\,{\rm where}\,\, M^{\O}=\bigoplus_{\gm\in \O}M^{\gm}.            
\end{equation}

Similarly, a weight $\mI_n$-module $M$ is a direct sum of its weight submodules $M_{\O}$,

\begin{equation}\label{MMO1}
M=\bigoplus_{\O\in\cM_n/G}M_{\O} \;\;{\rm where}\;\; M_{\O}=\bigoplus_{\gm\in \O}M_{\gm}.    
\end{equation}

For each orbit $\O\in \cM_n/G$, let $\WO$ (resp., $\GWO$) be the subcategory of weight (resp., generalized weight)
$\mI_n$-modules $M$ with $\Supp(M)\subseteq \O$. By (\ref{MMO}) and (\ref{MMO1}),

\begin{equation}\label{MMO2}
 \WIn=\bigoplus_{\O\in\cM_n/G}\WO \;\; {\rm and}  \;\; \GWIn=\bigoplus_{\O\in\cM_n/G}\GWO, 
\end{equation}

\noindent direct sum  of full subcategories of $\WIn$ and $\GWIn,$ respectively.

So, the problem of classification of indecomposable weight or generalized weight $\mI_n$-modules is reduced to the case when
the support of a module belongs to a single orbit.

Let $0\ra N\ra M\ra L\ra 0$ be a short exact sequence of $\mI_n$-modules. Then $M$ is a generalized weight module iff so are the modules $N$ and $L$, and in this case,

\begin{equation}\label{SMNL}
\Supp (M)= \Supp (N)\cup \Supp (L).
\end{equation}

{\bf The simple weight $\mI_n$-module $P_n$}.  By the  definition, the algebra $\mI_n$ 
is a subalgebra of the algebra $\End_K(P_n)$ of all 
$K$-endomorphism of the vector space $P_n$. 
So, the polynomial algebra $P_n$ is a (left)
$\mI_n$-module. Since $A_n\subset \mI_n$,
the $A_n$-module $P_n$ is a simple faithfull $A_n$-module.
Hence, the $\mI_n$-module $P_n$ is also simple and faithfull.
The action of the elements $x_i, \der_i, H_i$
and $\int_i$ on $P_n$ is given  by the rule: 
For all elements $p\in P_n$, 
$$\der_i p=\frac{\der p}{\der x_i},\;\; 
H_i p=\frac{\der}{\der x_i}(x_i p),\;\; 
\int_i p= \int_0^{x_i} pdx_i\;\; {\rm and}\;\; x_ip=x_i\cdot p 
\;\;({\rm multiplication\; by}\;\; x_i).$$

For all $\alpha= (\alpha_1, \ldots , \alpha_n)\in \N^{n}$
and $i=1, \ldots, n$, $H_ix^{\alpha}=(\alpha_i +1)x^{\alpha}$
where $x^{\alpha}=x_1^{\alpha_1}\cdots x_n^{\alpha_n}.$
Therefore, the $\mI_n$-module $P_n$ is a weight module,
$\Supp(P_n)=\N_+^{n}$. The polynomial algebra 
$P_n=K[x_1]\otimes \cdots \otimes K[x_n]$ is a tensor product
of polynomial algebras $K[x_i]$.
Futhermore, it is a tensor product of simple weight 
$\mI_1(i)$-modules $K[x_i]$ where 
$\mI_1(i)=K \langle \der_i , H_i, \int_i\rangle.$\\

{\bf The indecomposable $\mI_1$-modules $M(s,\lambda)$}, \cite{BBF-1}.
For $\lambda\in K$ and  a natural number $s\geq 1$, consider the $B_1$-module
\begin{equation}
 M(s,\lambda):= B_1\t_{K[H]} K[H)/(H-\lambda)^{s}. \label{FMM2}
\end{equation}

\noindent Clearly,
\begin{equation}
 M(s,\lambda)\simeq B_1/B_1 (H-\lambda)^{s}\simeq\mI_1/(F+\mI_1 (H-\lambda)^{s}). \label{FMM3}
\end{equation}

\noindent The $\mI_1$-module/$B_1$-module $M(s,\lambda)$ is a generalized weight module with $\Supp M(s,\lambda)=\lambda+\Z$,
\begin{equation}
 M(s,\lambda)=\bigoplus_{i\in \Z}M(s,\lambda)^{\lambda +i} \;\;{\rm and}\;\; \dim \, M(s,\lambda)^{\lambda+i}=s\;\; {\rm for\; all}\;\; i\in \Z . \label{FMM4}
\end{equation}

\begin{lemma}\label{a12Apr17}
 \begin{enumerate}
  \item Each simple generalized weight $\mI_1$-module is a simple weight $\mI_1$-module, and vice versa.
  
  \item Each simple generalized weight $\mI_1$-module is isomorphic to one of the modules: $K[x]$ or $M(1, \lambda),\;
  \lambda\in \Lambda$ (where $\Lambda$ is any fixed subset of $K$ such that the map  $\Lambda\to K/\Z,$ $\lambda\mapsto \lambda+\Z$
  is a bijection), they are pairwise non-isomorphic.
  
  \item $\End_{\mI_1}(M)\simeq K$ for all simple  weight $\mI_1$-modules.
 \end{enumerate}

\end{lemma}

{\it Proof}.
 1 and 2. Statements 1 and 2 follow from  \cite[Theorem 2.5]{BBF-1}.

3. $K \id_{K[x]}\subseteq \End_{\mI_1}(K[x])\subseteq \End_{A_1}(K[x])\simeq \Ker_{K[x]}(\der \cdot )\simeq K$
where $\der\cdot:K[x]\to K[x], p\mapsto \frac{dp}{dx}$ and
$A_1=K \langle x , \der\,|\,\der x-x \der=1 \rangle$ is the first Weyl algebra.
Therefore, $\End_{\mI_1}(K[x])\simeq K$.
Similarly, 
$$\End_{\mI_1}(M(1,\lambda))\simeq \End_{\mI_1}(\mI_1/(\mI_1(H-\lambda)))\simeq \Ker_{M(1, \lambda)}((H-\lambda) \cdot )\simeq K.\;\;\;\Box$$


\begin{lemma}\label{b12Apr17}

\begin{enumerate}
\item \cite[Proposition 6.1.(1)]{algintdif}  The $\mI_1$-module $K[x]$ is isomorphic to $\mI_1/\mI_1\der$. 
\item \cite[Eq. (12)]{BBF-1} For all elements $\l \in K\setminus\N_{+}$,  the $\mI_1$-module $M(1,\l)$ is isomorphic to the $\mI_1$-module $\mI_1/\mI_1(H-\l)$.
\end{enumerate}
\end{lemma}

For all $n\in \N_{+}$, 
$F+\mI_1(H_1-n)=E_{*,n-1}\oplus\mI_1(H_1-n)$,
where 
$E_{*,n-1}:=\oplus_{i\geq 0}Ke_{i,n-1}\simeq$ $_{\mI_1}K[x_i].$
Hence, $M(1,n)=\mI_1/(E_{*,n-1}\oplus\mI_1(H_1-n))$
and there is a short split sequence of $\mI_1$-modules

\begin{equation}\label{KxEI}
0\to K[x_1]\simeq E_{*,n-1}\to \mI_1/\mI_1(H_1-n)\to M(1,n)\to 0.
\end{equation}

\noindent In fact, it splits (Theorem \ref{B11Apr17}).

Let $A$ be an algebra and $\Md{A}$ be the category of left $A$-modules.
A subcategory $\CC$ of $\Md{A}$ is called a {\em semisimple} category if every 
module of $\CC$ is a direct sum of its simple modules in $\CC$.
The category $\GWIn$ of generalized weight $\mI_n$-modules is a subcategory of the category $\Md{\mI_n}$ of all left
$\mI_n$-modules.
The category $\WIn$ of weight $\mI_n$-modules is a subcategory of $\GWIn$.
For an $A$-module $M$, we denote by $\End_A(M)$ its algebra of endomorphisms.\\

{\bf Classification of simple weight $\mI_n$-modules.} We denote by $\widehat{\mI}_n({\rm weight})$ (resp., $\widehat{\mI}_n({\rm gen.\; weight})$)
the set of isomorphism classes of simple weight (resp., generalized weight) $\mI_n$-modules.
The next theorem classifies (up to isomorphism) all the simple weight $\mI_n$-modules.

\begin{theorem}\label{A11Apr17}
\begin{enumerate}
 \item $\widehat{\mI}_n(gen.\; weight)=\widehat{\mI}_n(weight)=\widehat{\mI}_1(weight)^{\otimes n}$, i.e.,
 any simple generalized weight $\mI_n$-module is a simple weight $\mI_n$-module, and vice versa;
 any simple weight $\mI_n$-module $M$ is isomorphic to the tensor product $M_1\otimes \cdots \otimes M_n$
 of simple weight $\mI_1$-modules and two such modules are isomorphic over $\mI_n$,
 $M_1\otimes \cdots \otimes M_n\simeq M_1^{\prime}\otimes \cdots \otimes M_n^{\prime},$
 iff for each $i=1, \ldots, n$, the $\mI_1$-modules $M_i$ and $M_i^{\prime}$ are isomorphic.
 
 \item For each simple weight $\mI_n$-module $M= \bigotimes_{i=1}^{n}M_i$, $\Supp(M)=\prod_{i=1}^{n}\Supp(M_i).$
\end{enumerate}

\end{theorem}

{\it Proof}. To prove the theorem we use induction on $n$.
The case $n=1$ is true, see  \cite[Theorem 2.5]{BBF-1}.

Suppose that $n>1$ and that theorem is true for all $n^{\prime}<n$.
Let $M$ be a simple generalized weight $\mI_n$-module.
By (\ref{MMO}), $\Supp(M)\subseteq \O$ for some orbit $\O\in K^{n}/\Z^{n}$.
Since $\mI_n=\mI_1\otimes \mI_{n-1}$, the $\mI_1$-module $M$ is a
generalized weight $\mI_1$-module with $\Supp_{\mI_1}(M)\subseteq \lambda_1+\Z$
where for some $\lambda_1\in K$ such that $\lambda=(\lambda_1, \ldots, \lambda_n)\in \Supp_{\mI_n}(M)$ for some
$\lambda_2, \ldots, \lambda_n\in K$.
Fix a nonzero vector $v\in M^{\lambda}$ such that $(H_1-\lambda_1)v=0$.
The $\mI_1$-submodule $\mI_1 v$ of $M$ is an epimorphic image of the $\mI_1$-module 
$N=\mI_1/\mI_1(H_1-\lambda_1)$.
By \cite{algintdifline}, the $\mI_1$-module $N$ has finite length. 
Hence, so is the  $\mI_1$-module $\mI_1 v$. 
Changing the element $v$, if
necessary, we can assume that the $\mI_1$-module $M_1=\mI_1 v$ is a simple weight $\mI_1$-module.
The $\mI_n$-module $M$ is a simple module. 
The $\mI_n$-module homomorphism
$$N\t \mI_{n-1}=\mI_1/\mI_1(H_1-\lambda_1)\t \mI_{n-1}\longrightarrow M,\;\;
\overline{1}\otimes 1\mapsto v\;\; ({\rm where }\;\; \overline{1}=1+\mI_1(H_1-\lambda_1))$$

\noindent is an epimorphism. 
By Lemma \ref{a12Apr17}.(3), $\End_{\mI_1}(M)\simeq K$.
By \cite{Bav-Schurian}, $M\simeq N\t M^{\prime}$ for some simple $\mI_{n-1}$-module  $M^{\prime}$.
The $\mI_n$-module $M$ is a generalized weight $\mI_n$-module.
Hence,  the $\mI_{n-1}$-module $M^{\prime}$ is a generalized weight $\mI_{n-1}$-module.
Now, the result follows by induction on $n$.
$\Box$

\

{\bf The category $\WIn$ of weight $\mI_n$-modules is a semisimple category.}
Let $n=1$. For each orbit $\O\in K/\Z$, we fix an element
$\lambda_{\O}\in K$ such that $\lambda_{\O}\in \O$.
In particular, $\lambda_{\O}+\Z=\O$. For $\O=\Z$ let 
$\lambda_{\Z}:=0$. Let $n>1$. for each orbit 
$\O=\O_1\times\cdots \times \O_n\in K^{n}/\Z^{n}$,
let $\lambda_{\O}:=(\lambda_{\O_1}, \ldots, \lambda_{\O_n})
\in K^{n}$. In particular $\lambda_{\O}+\Z^{n}=\O$.
The map $\O\mapsto \lambda_{\O}$ is a bijection, by definition.
An orbit $\O=\O_1\times\cdots\times\O_n\in K^{n}/\Z^{n}$
is a direct product of orbits $\O_i\in K/\Z$. For the orbit
$\O$, let 
$$\mD_{\O}:= \{i\in\{1,\ldots, n\}\,|\,\O_i=\Z\}.$$
Then $\{1,\ldots,n\}\setminus \mD_{\O}=\{j\in\{1,\ldots,n\}\,|\,\O_j\neq \Z\}.$
Let $\CD_{\O}$ be any subset of $\mD_{\O}$
(eg.,  $\CD_{\O}=\emptyset)$.
Then

\begin{equation}\label{MOMO1}
 \{1, \ldots, n\}=\CD_{\O}\sqcup \CN_{\O}
\end{equation}
\noindent is a disjoint union where $\CN_{\O}=\{1, \ldots, n\}\setminus \CD_{\O}$.
For each pair $(\O, \CD_{\O})$, let us define the $\mI_n$-module
\begin{equation}\label{MOMO}
M(\CD_{\O}):=\bigotimes_{i=1}^{n}M(\CD_{\O_i})\;\;
{\rm where} \;\; M(\CD_{\O_i})=\begin{cases}
K[x_i]& \text{if } i\in \CD_{\O},\\
M(1,\l_{\O_i})& \text{if } i\not\in \CD_{\O}.\\
\end{cases}
\end{equation}



\noindent For each choice of the set $\CD_{\O}$, $\O=\Od\times\Ond$ where $\Od:=\prod_{i\in\CD_{\O}}\O_i\simeq \N_{+}^{|\CD_{\O}|}$ is called
the {\em degenerate part} of the pair $(\O, \CD_{\O})$ and $\Ond:=\prod_{j\in\CN_{\O}}\O_j$ is called
the {\em non-degenerate part} of the pair $(\O, \CD_{\O})$.
The elements of the set $\CD_{\O}$ (resp., $\CN_{\O}$) are called the {\em degenerate} 
(resp., {\em non-degenerate}) {\em indices} with respect to the pair $(\O, \CD_{\O})$, or, simply,
the $\CD_{\O}$-{\em degenerate} (resp., $\CD_{\O}$-{\em non-degenerate}) {\em indices}.

For each subset $I$ of $\{1,\ldots,n\}$, let $\ga_n(I)$ be the ideal of $\mI_n$ generated by the ideals $F(i)$ of $\mI_1(i)$ where $i\in I$.
If $I=\emptyset$ we set $\ga_n(\emptyset)=0$.
Clearly, $\ga_n(I)=\sum_{i\in I}\gp_{i}.$
If $I=\{1, \ldots, n\}$ then $\ga_n(I)=\ga_n$ is the maximal ideal
of $\mI_n$.
The factor algebra
$$\mI_n(\CD_{\O}):=\mI_n/\ga_n(\CD_{\O})\simeq B(\Od)\otimes \mI_l(\Ond)$$
\noindent is a tensor product of algebras where
$$B(\Od):=\bigotimes_{i\in{\CD_{\O}}}B_1(i),\;\; \mI_l(\Ond):=\bigotimes_{i\in{\CN_{\O}}}\mI_1(i) \;\;{\rm and}\;\; l=|\CN_{\O}|.$$

For an algebra $A$ and an $A$-module $M$, $\ann_{A}(M):=\{a\in A\,|\,aM=0\}$ is
the {\em annihilator} of the $A$-module $M$. 
The $A$-module $M$ is called a {\em faithfull} $A$-module if $\ann_A(M)=0$.

\begin{corollary}\label{a15Apr17}
 Let $\O\in K^{n}/\Z^{n}$. Then the set $\GWO^{\widehat{}}$ of isomorphisms classes of  simple (generalized) weight $\mI_n$-modules
 in the category $\GWO$ is equal to $\{M(\CD_{\O})\,|\, \CD_{\O}\subseteq \mD_{\O}\}$ and  the number of elements in this
 set is $2^{|\mD_{\O}|}$.
 \begin{enumerate}
  \item The $\mI_n$-module $M(\CD_{\O})$ is faithfull iff $\CD_{\O}=\{1,\ldots,n\}$.
  \item If $\CD_{\O}\neq\{1,\ldots,n\}$ then $\ann_{\mI_n}(M(\CD_{\O}))=\ga_n(\CN_{\O})$.
 \item The elements in $\GWO^{\widehat{}}$ are uniquely determined by their annihilators, i.e., the map $M(\CD_{\O})\mapsto \ann_{\mI_n}(M(\CD_{\O}))$ is a bijection.
 
 \end{enumerate}

\end{corollary}

{\it Proof}. Since $\mI_n=\t_{i=1}^{n}\mI_1(i)$ and $M(\CD_{\O})=\t_{i=1}^{n}M(\CD_{\O_i})$,
we have $\ann_{\mI_n}(M(\CD_{\O}))=\ga_n(\CN_{\O})$ and statements 1 and 2 follow. Then statement 2 implies statement 3, and the rest follows. \;\;$\Box$

 Clearly,
\begin{equation}
 S(\CD_{\O}):=\Supp(M(\CD_{\O}))=\N_{+}^{|\CD_{\O}|}\times \prod_{i\in C\CD_{\O}} \O_i 
 \label{SDCO}
\end{equation}
\noindent where $\O=\prod_{i=1}^{n}\O_i$ and 
$C\CD_{\O}=\{1,\ldots, n\}\setminus \CD_{\O}$
is the complement of the set $\CD_{\O}$ in $\{1,\ldots, n\}$.

\begin{theorem}\label{B11Apr17}
Every weight $\mI_n$-module is a direct sum of simple weight 
$\mI_n$-modules. In particular,
the category $\WIn$ of weight $\mI_n$-modules is a semisimple category.
\end{theorem}

{\it Proof}.  In view of (\ref{MMO2}), it suffices to show that
for each orbit $\O\in K^{n}/\Z^{n}$ and  for any two weight simple 
$\mI_n$-modules $M$ and $\overline{M}$  with
$\Supp(M), \Supp(\overline{M})\subseteq \O$ the short exact sequence 

$$ 0\to M\to M^{\prime}\stackrel{f_{\lambda}}{\to} \overline{M}\to 0  $$

\noindent splits where $M^{\prime}$ is a weight $\mI_n$-module.
Notice that $\Supp(M^{\prime})=\Supp(M)\bigcup \Supp(\overline{M})\subseteq \O$.
By  Corollary \ref{a15Apr17}, $M=M(\CD_{\O})$ and $\overline{M}=M(\overline{\CD}_{\O})$ for some subsets $\CD_{\O}$ and 
$\overline{\CD}_{\O}$ of the set $\mD_{\O}$.
 
 For each element $\l = (\l_1, \ldots ,\l_n)\in \Supp (\bM )\subseteq\OO$, 
 we have the short exact sequence of $K_{\gm_\l}$-modules where $K_{\gm_\l}:= D_n/\gm_\l$ and $\gm_\l = (H_1-\l_1, \ldots ,  H_n-\l_n)$, 
$$ 0 \ra M_\l\ra M_{\l}' \stackrel{f}{\ra}\bM_\l \ra 0.$$ 
Notice that $\dim_K(M_\l) = \dim_K(\bM_\l ) =1$. 
Since the algebra $K_{\gm_\l}$ is a field, the short exact sequence above splits. Let $\overline{e}=\overline{e}_\l$ 
be a nonzero element of $\bM_\l$ and $e=e_\l$ be an element of $M_\l'$ such that $f_\l (e) = \overline{e}$.
To finish the proof it suffices to show that the $\mI_n$-submodule of $M'$, $N:=\mI_ne$, is a simple $\mI_n$-module: 
Indeed, in this case, $M\cap N=0$ since $f(M)=0$, $f(N)=\bM$ and the $\mI_n$-modules $M$ and $N$ are simple. 
Then the $\mI_n$-module $M'$ of length 2 contains the submodule $M\oplus N$ of length 2, and so $M'=M\oplus N$, as required.
There are two cases to consider: either $\overline{\CD}_\OO = \emptyset$ or $\overline{\CD}_\OO  \neq  \emptyset$.

(i) $\overline{\CD}_\OO = \emptyset$: In this case, $\bM =\bigotimes_{i=1}^n\bM_i$ where $\bM_i=\mI_1(i)/ \mI_1(i)(H_i-\l_i)$ and $\l_i\in K\backslash \Z$. 
For each number $i=1, \ldots , n$, we have the $\mI_1(i)$-module epimorphism 
$$\bM_i=\mI_1(i)/ \mI_1(i)(H_i-\l_i)\ra \mI_1(i)e, \;\;  1+\mI_1(i)(H_i-\l_i)\mapsto e,$$ which is
necessarily an isomorphism since the $\mI_1(i)$-module $\bM_i$ is simple.
Since $\mI_n = \mI_1(1)\t \cdots \t \mI_1(n)$, we have the $\mI_n$-epimorphism 
$$\bM=\bigotimes_{i=1}^n\bM_i=\bigotimes_{i=1}^{n}\mI_1(i)/ \mI_1(i)(H_i-\l_i)\ra \mI_ne=N$$
\noindent which is necessarily an isomorphism since the  $\mI_n$-module $\bM$ is simple. 

(ii) $\overline{\CD}_\OO  \neq  \emptyset$: Without loss of generality we can assume that $\overline{\CD}_\OO  = \{ 1, \ldots , m\}$ for some natural number $m$ such that $1\leq m \leq n$. Then 
$$ \bM_i=\begin{cases}
K[x_i]& \text{if }i=1, \ldots , m,\\
M(1,\l_i)& \text{otherwise}.\\
\end{cases}
$$
Notice that $\Supp_{\mI_1(i)}(K[x_i])=\{ 1,2,\ldots \}$. Choose $\l=(\l_1,\ldots, \l_n)\in\Supp(\overline{M})$ such that
$\l_i=1$ for all $i=1,\ldots ,m$. Then 
$$ H_ie=\begin{cases}
e& \text{if }i=1, \ldots , m,\\
\l_ie& \text{otherwise}.\\
\end{cases}
$$
For each $i=1,\ldots ,m$, 
$$0=(H_i-1)e=(\der_ix_i-1)e=(x_i\der_i+[\der_i, x_i]-1)e=(x_i\der_i+1-1)e=x_i\der_ie.$$
Since $\bM=\bigotimes_{j=1}^n\bM_j$ and $\bM_i=K[x_i]$ for $i=1,\ldots , m$, the map $x_i\cdot : \bM\ra \bM$, $ m\mapsto x_i m$ is an injection. 
Therefore, $\der_ie=0$ for $i=1,\ldots , m$. So, by Lemma \ref{a12Apr17}.(1), we have an $\mI_1(i)$-epimorphism 
$$ K[x_i]=\mI_1(i)/\mI_1(i)\der_i\ra \mI_1(i) e, \;\; 1+\mI_1(i)\der_i\mapsto e,$$ which 
is necessarily an isomorphism since the $\mI_1(i)$-module $K[x_i]$ is simple. Using the same argument as in the case (i),
we see that
for all $j>m$, $\mI_1(j)e\simeq M(1,\l_j)$ and the $\mI_n$-epimorphism $\bM=\bigotimes_{i=1}^n\bM_i\ra \mI_ne=N$ is an isomorphism since the $\mI_n$-module $M$ is simple.
The proof of the theorem is complete. $\Box $

By Theorem \ref{B11Apr17}, each weight $\mI_n$-module
$M$ is a unique direct sum 

\begin{equation}
 M=\bigoplus_{\O\in K^{n}/\Z^{n}}\bigoplus_{\CD_{\O}\subseteq
 \mD_{\O}} M(\CD_{\O})^{\mu(\CD_{\O})} \label{MWD}
\end{equation}
\noindent where $\mu(\CD_{\O})$ is the multiplicity of
$M(\CD_{\O})$ (which can be any set).


\section{Explicit description of indecomposable generalized weight  $\mI_n$-modules}

In this section, an explicit description of indecomposable generalized weight  $\mI_n$-modules is obtained (Theorem \ref{14Apr17}). 
One of the key steps is to show that each category $\GWO$ is
a direct sum of its subcategories $\GWDO$ that are generated
by the single simple weight $\mI_n$-modules $M(\CD_{\O})$, see (\ref{GWDD}).
Using (\ref{GWDD}) and some results about representations of  Artinian rings,  a criterion is given for the category $\GWO$
and its subcategories $\GWDO$ to be of finite representation type, tame or wild. Explicit classes of indecomposable modules in  $\GWO$ are considered.

Let $A$ be an algebra, $\gm$ be a {\em co-finite ideal} of $A$ (i.e., $\dim_K(A/\gm)<\infty)$.
 An $A$-module $M$ is called a {\em locally finite module} if for each $m\in M$, $\dim_K(Am)<\infty$.
 An $A$-module $M$ is called an {\em $\gm$-locally finite module} if for each $m\in M$, $\dim_K(Am)<\infty$
 and $\ann_{A}(Am)\subseteq \gm^{i}$ for some $i\geq 1$ (if, in addition, $A$ is a commutative algebra then the last condition
 is equivalent to the condition that $\gm^{i}m=0$ for some $i\geq 1$).
We denote by $\LFm(A)$ the category of all $\gm$-locally finite $A$-modules.
The category $\LFm(A)$ is closed under arbitrary direct sums, submodules and factor modules.
An $\mI_n$-module $M\in \GWO$ (resp., $M\in \WO$) is called {\em equidimensional} 
if $\dim_K(M^{\gm})=\dim_K(M^{\gn})$ (resp., $\dim_K(M_{\gm})=\dim_K(M_{\gn})$)
for all $\gm, \gn\in \Supp(M)$. If the common value of all $\dim_K(M^{\gm})$ (resp., $\dim_K(M_{\gm})$) is $d$,
we say that $M$ is $d\mbox{-}${\em equidimensional}.
Let $\gm$ be a maximal ideal of the polynomial algebra $D_n$ and $\CI(D_n, \gm)$ is the set of all ideals
$I$ of $D_n$ such that $\gm\supseteq I\supseteq \gm^{i}$ for some $i\geq 1$.
For all ideals $I\in\CI(D_n, \gm)$, the factor algebra $D_n/I$ is a local, finite dimensional, commutative algebra
with maximal ideal $\gm/I$.
\\

{\bf The category $\GWO$ is a direct sum of subcategories $\GWDO$.} Let 
$\GWDO$ be the full subcategory of $\GWO$ generated by the simple
weight $\mI_n$-module $M(\CD_{\O})$.
There are precisely $2^{|\mD_{\O}|}$ such subcategories in the 
categories $\GWO$.
They are key objects in the description of all indecomposable generalized weight modules $M$ with $\Supp(M)\subseteq \O$ 
since

\begin{equation}\label{GWDD}
 \GWO =\bigoplus_{\CD_{\O}\subseteq \mD_{\O}} \GWDO ,
\end{equation}

\noindent Theorem \ref{15May19}.

Let $R$ be a ring, $M$ be a nonzero $R$-module and 
$\gp=\ann_{R}(M)$. The $R$-module $M$ is called a 
{\em prime} $R$-module (or a $\gp$-{\em prime} $R$-module)
if $\gp$ is a prime ideal of $R$ and $\ann_R(N)=\gp$ 
for all nonzero $R$-submodules $N$ of $M$.
The $R$-module $M$ is called an {\em absolutely prime}
$R$-module (or an {\em absolutely} $\gp$-{\em prime} $R$-module)
if $\gp$ is a prime ideal of $R$ and $\ann_R(N)=\gp$ for all
nonzero subfactors of $M$, i.e., $N=M_2/M_1$ for some submodules $M_1$ and $M_2$ of $M$ such that $0\subseteq M_1\subset M_2\subseteq M$.

Lemma \ref{a15May19} shows that all nonzero modules in each category $\GWDO$ are absolutely prime $\mI_n$-modules 
(Corollary \ref{a21May19}).
Lemma \ref{a15May19} is one of the key steps in proving that the equality (\ref{GWDD}) holds.

\begin{lemma}\label{a15May19}
Let $0\neq M\in \GWDO$. Then

\begin{enumerate}
\item $\Supp(M)=\Supp(M(\CD_{\O}))$.
\item $\ann_{\mI_n}(M)=\ga(\CD_{\O})$ where
$\ga(\CD_{\O})=\ann_{\mI_n}(M(\CD_{\O}))=\sum_{i\in C\CD_{\O}}
\gp_i\in \Spec(\mI_n)$ where $C\CD_{\O}=\{1,\ldots, n\}\setminus \CD_{\O}$.
So, all nonzero modules in the category $\GWDO$ are absolutely 
$\ga(\CD_{\O})$-prime $\mI_n$-modules.
\end{enumerate}
\end{lemma}

{\it Proof}. 
 1. By the definition, the category of generalized weight
$\mI_n$-modules $\GWDO$ is generated by the simple weight
$\mI_n$-module $M(\CD_{\O})$, and statement 1 follows, by
(\ref{SMNL}). 

2. The ideal 
$$\ga(\CD_{\O})=\ann_{\mI_n}(M(\CD_{\O}))=
\sum_{i\in C\CD_{\O}}\gp_i$$ 
is a prime ideal of $\mI_n$.
Let $N=\ann_{M}(\ga)=\{m\in M\,|\,\ga m=0\}$ where 
$\ga=\ga(\CD_{\O})$. We have to show that $M=N$.
Suppose not, we seek a contradiction. 
Then $S:=\soc_{\mI_n}(M/N)\neq 0$ since 
$0\neq M/N\in\GWDO$. Now, $S=L/N$ for some submodule $L$ of $M$ such that $N\varsubsetneq L$. 
Recall that $\ga^{2}=\ga$. 
So, $0\neq \ga L=\ga^{2}L\subseteq \ga N=0$, 
a contradiction.\;\;$\Box $

\begin{theorem}\label{15May19}
For all $\O\in K^{n}/\Z^{n}$ and all $\CD_{\O}\subseteq\mD_{\O}$,
the equality (\ref{GWDD}) holds.
\end{theorem}

{\it Proof}. {\sc Claim.} 
$\sum_{\CD_{\O}\subseteq\mD_{\O}} \GWDO=\bigoplus_{\CD_{\O}\subseteq\mD_{\O}} \GWDO$.

Suppose that the Claim does not hold, i.e., 
$M_1+ \cdots +M_s\neq M_1\oplus\cdots \oplus M_s$ for some 
modules $M_i\in \GWDOi$  such that the sets $\CD_{\O_1},\ldots , \CD_{\O_s}$
are distinct.
Then $s\geq 2$.
By Corollary \ref{a15Apr17}, the prime ideals
$$\ga_1=\ann_{\mI_1}(M_1), \ldots, \ga_s=\ann_{\mI_1}(M_s)$$
\noindent are {\em distinct}. 
Up to order, we may assume that the ideal $\ga_1$ a minimal element (with respect to $\subseteq$) of the set 
$\{\ga_1,\ldots, \ga_s\}$.
By replacing the module $M_1$, by a (possibly) smaller 
nonzero submodule, say $M_1^{\prime}\subseteq M_1$, 
we may assume that $M_1^{\prime}\subseteq M_2+\cdots +M_s$.
By Lemma \ref{a15May19}.(2), $\ann_{\mI_1}(M_1^{\prime})=\ga_1$.
Let $\ga=\ga_2\cdots\ga_s$. 
Then $\ga M_1^{\prime}\subseteq \ga(M_2+\cdots+M_s)=0$,
and so 
$\ga_1\cdots \ga_s\subseteq \ann_{\mI_1}(M_1^{\prime})=\ga_1.$
Recall that the ideal $\ga_1$ is a prime ideal (Lemma \ref{a15May19}.(2)), hence $\ga_i\subseteq\ga_1$ 
for some $i$ such that $2\leq i\leq s$.
This contradicts the minimality of the ideal $\ga_1$.
Now, the theorem follows from Proposition \ref{A15May19}.
\;\;$\Box $

\begin{proposition}\label{A15May19}
Let $\O\in K^{n}/\Z^{n}$. 
Suppose that $\CD_{\O}$ and $\CD_{\O}^{\prime}$ are distinct
subsets of the set $\mD_{\O}$. 
Let $M=M(\CD_{\O})$ and $M^{\prime}=M(\CD_{\O}^{\prime}).$
Then $\Ext_{\mI_n}^{1}(M^{\prime}, M)=0.$
\end{proposition}

{\it Proof}. By Lemma \ref{a15May19}.(2), the ideals
$\ga=\ann_{\mI_n}(M)$ and $\ga^{\prime}=\ann_{\mI_n}(M^{\prime})$
are {\em distinct} prime ideals of the algebra $\mI_n$.
Therefore, either $\ga\nsubseteq \ga^{\prime}$ or
otherwise $\ga\subset \ga^{\prime}$ (a proper inclusion since 
$\ga\neq \ga^{\prime}$).
Let $0\to M\to N\to M^{\prime}\to 0$ be a short exact sequence
of $\mI_n$-modules. 
To finish the proof it suffices to show that the short exact sequence splits.

 (i) {\em Suppose that} $\ga\nsubseteq\ga^{\prime}$:
 Then $\ga M^{\prime}\neq 0$ (since otherwise $\ga M^{\prime}=0$,
 and so $\ga\subseteq \ann_{\mI_n}(M^{\prime})=\ga^{\prime}$,
 a contradiction).
 In particular, $\ga M^{\prime}=M^{\prime}$ since 
 the $\mI_n$-module $M^{\prime}$ is simple.
 Now, 
 $$N\supseteq \ga N=\ga(N/M)=\ga M^{\prime}=M^{\prime}.$$
 Therefore, the $\mI_n$-submodule $\ga N$ of $N$ 
 is isomorphic to the simple $\mI_n$-module $M^{\prime}$.
 Hence, $N\supseteq M+\ga N=M\oplus\ga N\simeq M\oplus M^{\prime}$,
 i.e., $N=M\oplus\ga N\simeq M\oplus M^{\prime}$,
 i.e., the short exact sequence splits.
 
 (ii) {\em Suppose that} $\ga\subset \ga^{\prime}$:
 Let $\CD=\CD_{\O}$ and $\CD^{\prime}=\CD_{\O}^{\prime}$.
 By Lemma \ref{a15May19}.(2),
 $\ga=\sum_{i\in C\CD}\gp_i\subset \ga^{\prime}=\sum_{j\in C\CD^{\prime}}\gp_j$, and so $C\CD\subset C\CD^{\prime}$ or,
 equivalently, $\CD\supset \CD^{\prime}$.
 Up to order, let $\mD_{\O}=\{1,\ldots, m\}$,
 $\CD=\{1,\ldots,l\}\coprod \CD^{\prime}$ and 
 $\CD^{\prime}=\{l+1,\ldots, k\}$ provided 
 $\CD^{\prime}\neq \emptyset .$
 Notice that $1\leq l\leq m$ and $k\leq m$, 
 $\O=\Z^{m}\times (\lambda_{m+1}+\Z)\times\cdots \times(\lambda_n+\Z)$ where $\lambda_i\notin \Z$ for all $i$ such that $m+1\leq i\leq n.$

Clearly, 
\begin{eqnarray*}
\Supp (M) &= &\N_+^k\times \Z^{m-k}\times \prod_{i=m+1}^n(\l_i+\Z ),\\
\Supp (M') &=& \Z^l\times \N_+^{k-l}\times \Z^{m-k}\times \prod_{i=m+1}^n(\l_i+\Z )
\end{eqnarray*}
since 
\begin{eqnarray*}
M &=& K[x_1, \ldots , x_k]\otimes \bigotimes_{j=k+1}^m \mI_1(j)/\mI_1(j)H_j\otimes \bigotimes_{i=m+1}^n\mI_1(i)/\mI_1(i)(H_i-\l_i),\\
M' &=&\bigotimes_{s=1}^l\mI_1(s)/\mI_1(s)H_s\otimes  K[x_{l+1}, \ldots , x_k]\otimes \bigotimes_{j=k+1}^m \mI_1(j)/\mI_1(j)H_j\otimes \bigotimes_{i=m+1}^n\mI_1(i)/\mI_1(i)(H_i-\l_i).
\end{eqnarray*}

Let $\l = (1,\ldots , 1, \l_{m+1},\ldots , \l_n)$. Then $\dim_K(M^\l )=\dim_K(M'^\l )=1$ and $\dim_K(N^\l )=2$ since there is a short exact sequence $0\ra M^\l \ra N^\l \ra M'^\l\ra 0$ of $K$-modules. Fix an element $v\in N^\l \setminus M^\l$. Then $v\neq 0$. Let $\th :=(1-e)v$ where $e=e_{00}(1)=1-\int_1\der_1$. Then $\th \equiv v \mod M$ since $eM'=0$. In particular, $\th \neq 0$.

(iii) $\int_1\der_1\th =\th$: Notice that $\int_1\der_1=1-e$ is an idempotent and the result follows:
$$ \int_1\der_1\th =(1-e)\th = (1-e)(1-e)v =(1-e)v  =\th.$$

(iv) $\der_1\th \neq 0$: Since $\th \neq 0$, the statement (iv) follows from the statement (iii).

(v) $H_1\der_1\th =0$: Since $(H_1-1)M'^\l =0$, we must have $ (H_1-1)N^\l \subseteq M^\l$. Now, 
$$ H_1\der_1\th = \der_1(H_1-1)\th \in \der_1(H_1-1)N^\l \subseteq \der_1M^\l=0.$$

The $\mI_1(1)$-module  $\mI_1(1)/\mI_1(1)H_1$ is a simple weight module. By the statement (v), the $\mI_1(1)$-submodule $L_1=\mI_1(1)\der_1\th$ of $N$ is isomorphic to 
the $\mI_1(1)$-module  $\mI_1(1)/\mI_1(1)H_1$. Recall that $\End_{\mI_1(1)}(L_1)=K$ and $\mI_n = \mI_1\t \mI_{n-1}$. By 
\cite{Bav-Schurian}, the $\mI_n$-submodule $\mI_n\der_1\th$ of $N$ is isomorphic to the tensor product $L_1\t L$ of the $\mI_1$-module $L_1$ and an $\mI_{n-1}$-module $L$.

(vi) $N=M\oplus L_1\t L$: Since the map $\der_1\cdot : L_1\t L\ra L_1\t L$, $u\mapsto \der_1u$ is a bijection (since the map $\der_1\cdot : L_1\ra L_1$, $w\mapsto \der_1w$ is so) and the map $\der_1\cdot : M\ra M$, $p\mapsto \der_1p$ has nonzero kernel (since $\der_1 M^\l =0$ and $M^\l \neq 0$), the simple $\mI_n$-module $M$ is not a submodule of $L_1\t L$. Hence, $M\cap L_1\t L =0$, and so $ M\oplus L_1\t L\subseteq N$, and the statement (vi) follows since the length of the $\mI_n$-module $N$ is 2.

Now, the proposition follows from the statement (vi). $\Box$

\begin{corollary}\label{a21May19}
 \begin{enumerate}
  \item Every module $M\in\GWO$ is a unique direct sum of
  absolutely prime generalized weight $\mI_n$-modules, 
  and this direct sum is $M=\bigoplus_{\CD_{\O}\subseteq \mD_{\O}}M_{\CD_{\O}}$ where $M_{\CD_{\O}}\in \GWDO.$
  \item Every generalized weight module is a unique sum of
  absolutely prime (generalized weight) $\mI_n$-modules.
 \end{enumerate}

 \end{corollary}

 {\it Proof}. 1. Statement 1 follows from Corollary \ref{a15Apr17} and Theorem \ref{15May19}.
 
 2. Statement 2 follows from Statement 1 and 
 (\ref{MMO2}). \;\;$\Box$
 
The next proposition shows that there are plenty of indecomposable generalized weight $\mI_n$-modules with  support
from a single orbit.

\begin{proposition}\label{A15Apr17}
Let $\O\in \cM_n/G$ and $\gm\in\O$.
\begin{enumerate}
 \item If $\mD_{\O}=\emptyset$ then $\{V(I):=B_n\t_{D_n}D_{n}/I\,|\,I\in\CI(D_n,\gm)\}\subseteq \GWDO$ 
 and each $\mI_n$-module $V(I)$ is an indecomposable, equidimensional, generalized weight $\mI_n$-module
 of length $d:=\dim_K(D_n/I)<\infty$ which is isomorphic to the $\mI_n$-module $B_n/B_nI=\bigoplus_{\alpha\in\Z^{n}}\der^{\alpha}D_n/I$
 (where $\der^{\alpha}=\der_1^{\alpha_1}\cdots\der_n^{\alpha_n}$) and $\dim_K(V(I)^{\gn})=d$ for all $\gn\in\O$, and $\ann_{\mI_n}(V(I))=\ga_n$.
 \item Suppose that $\mD_\O =\{ 1,\ldots , l\}$ (up to order) for some $l$ such that $1\leq l \leq n$. 
 \begin{enumerate}
 \item Suppose that $\CD_{\O}=\emptyset$. Let $\gm =(H_1,\ldots , H_l, H_{l+1}-\l_{l+1}, \ldots, H_n-\l_n)\in \Supp(M(\CD_\O) )$. Then  $\{V(I):=B_n\t_{D_n}D_{n}/I\,|\,I\in\CI(D_n,\gm)\}\subseteq \GWDO$ 
 and each $\mI_n$-module $V(I)$ is an indecomposable, equidimensional, generalized weight $\mI_n$-module
 of length $d=\dim_K(D_n/I)<\infty$ which is isomorphic to the $\mI_n$-module $B_n/B_nI=\bigoplus_{\alpha\in\Z^{n}}\der^{\alpha}D_n/I$ and $\dim_K(V(I)^{\gn})=d$ for all $\gn\in\Supp (V(I))=\Supp (M(\CD_\O ))$, and $\ann_{\mI_n}(V(I))=\ga_n$.
 \item  If $\CD_{\O}\neq\emptyset$ then $\CD_{\O}=\{1, \ldots, m\}$, up to order, for some $m$ such that
 $1\leq m\leq l$. Let $k=n-m$, $D_k=K[H_{m+1}, \ldots, H_n]$ and $B_k=\t_{i=m+1}^{n}B_1(i)$.
 Then 
 $$\{ V(I):=P_m\t(B_k\t_{D_k} D_k/I)\,|\,I\in\CI(D_k,\gm')\}\subseteq\GWDO$$
 \noindent where $\gm'=(H_{m+1},\ldots , H_l, H_{l+1}-\l_{l+1}, \ldots, H_n-\l_n)\in S(\CD_\O )$, and each
  $\mI_n$-module $V(I)$ is an indecomposable, equidimensional, generalized weight $\mI_n$-module
 of length $d:=\dim_K(D_k/I)<\infty$ which is isomorphic to the $\mI_n$-module $$P_m\t B_k/B_kI\simeq\bigoplus_{\alpha\in\Z^{k}}P_m\t \der^{\alpha}D_k/I$$
 (where $\der^{\alpha}=\der_{m+1}^{\alpha_1}\cdots\der_n^{\alpha_k}$) and $\dim_K(V(I)^{\gn})=d$ for all $\gn\in \Supp(V(I))=\Supp (M(\CD_\O ))$, and $\ann_{\mI_n}(V(I))=\ga_n(\CD_\O)$.
\end{enumerate}
\end{enumerate}
 
\end{proposition}

{\it Proof}. 1. Let $I\in\CI(D_n,\gm)$. Notice that $B_n=\bigoplus_{\alpha\in\Z^{n}}\der^{\alpha}D_n.$
Then 
$$V(I)=\bigoplus_{\alpha\in\Z^{n}}\der^{\alpha}D_n\t_{D_n}D_n/I\simeq \bigoplus_{\alpha\in\Z^{n}}\der^{\alpha}\t D_n/I=\bigoplus_{\alpha\in\Z^{n}}\der^{\alpha}D_n/I.$$
So, the $\mI_n$-module $V(I)$ is an equidimensional, generalized weight module with $\Supp(V(I))=\O$ and $\dim_K(V(I)^{\gn})=l=\dim_K(D_n/I)<\infty$.
By Theorem \ref{A11Apr17}, the  simple (generalized) weight $\mI_n$-module $M(\CD_\O)$ from $\GWO$ has support $\O$ and is $1$-equidimensional. Hence, 
$l_{\mI_n}(V(I))=\dim_K(D_n/I)=l$. 

It remains to show that the $\mI_n$-module $V(I)$ is an indecomposable. 
The functor 
$$B_n\t_{D_n} -:\Md{D_n}\rightarrow\Md{B_n},\;\;N\mapsto B_n\t_{D_n}N$$
\noindent is exact. The commutative algebra $D_n/I$ is a local, commutative, finite dimensional algebra with maximal ideal $\gm/I$.
Since $(D_n/I)/(\gm/I)\simeq D/\gm$ is a field, the $D_n$-module $D_n/I$ is indecomposable. Hence, so is the induced module $V(I)=B_n\t_{D_n}D_n/I$.
Clearly, $\ga_n\subseteq \ann_{\mI_n}(V(I))$. Since $\ga_n$ is a maximal ideal of $\mI_n$ and $V(I)\neq 0$, we must have $ \ann_{\mI_n}(V(I))=\ga_n$.

2(a). Repeat the arguments of statement 1.

2(b). The functor
$$P_m\t -: \Md{\mI_k}\rightarrow \Md{\mI_n}, \;\;L\mapsto P_m\t L$$
\noindent is an exact functor. Now, statement 2 follows from statement 2(a).\;\;$\Box$\\

{\bf Explicit description of modules in $\GWO$.}
In view of Theorem \ref{15May19}, Theorem \ref{14Apr17} below is an explicit description of generalized weight $\mI_n$-modules.

\begin{theorem}\label{14Apr17}
Let $\OO \in \cM_n/G$.
\begin{enumerate}
\item Suppose that $\CD_\OO= \emptyset$,   $\gm\in \OO$   if $\mD_\O =\emptyset$ and $\gm =(H_1,\ldots , H_l, H_{l+1}-\l_{l+1}, \ldots, H_n-\l_n)\in S(\CD_\O )$ if $\mD_\OO = \{ 1, \ldots , l\}\neq \emptyset$, up to order. Then the functor
$$ \GWDO \ra \LFm (D_n), \;\; M\mapsto M^\gm$$ is an equivalence of categories with the inverse $N\mapsto B_n\t_{D_n}N$, the induced functor. 
\item Suppose that $\CD_\OO\neq \emptyset$ and, up to order, $\CD_\OO = \{ 1, \ldots , m\}$, $\mD_\OO =\{ 1, \ldots , l\}$ for some $m$ such that $1\leq m \leq l\leq n$. Then $\mI_n=\mI_m\t \mI_{n-m}$ and 
$\OO = \Z^m\times \O'$ where $\O'=\Z^{l-m}\times\OO_1\times\cdots\times  \OO_{n-l}$ and $\OO_i\neq \Z$ for all $i=1,\ldots ,n-l$, and 
\begin{enumerate}
\item $\GWDO = P_m\t \GW (\mI_{n-m}, \CD_{\O'}'):=\{ P_m\t M\, | \, M\in \GW (\mI_{n-m},  \CD_{\O'}')\}$ with $ \CD_{\O'}'=\emptyset $.
\item Fix $\gm \in \OO$ such that $\gm =(H_1-1,\ldots , H_m-1, H_{m+1}, \ldots, H_l, H_{l+1}-\l_{l+1},\ldots , H_n-\l_n)$. 
Then $\gm' =(H_{m+1}, \ldots, H_l, H_{l+1}-\l_{l+1},\ldots , H_n-\l_n)\in\O'$ and  the functor 
$$ \GWDO \ra \LF_{\gm'} (D_{n-m}), \;\; P_m\t M\mapsto M^{\gm'}$$
is an equivalence of categories with the inverse $N\mapsto P_m\t(B_{n-m}\t_{D_{n-m}}N)$ where $D_{n-m}=K[H_{m+1},\ldots , H_n]$ and $D_0:=K$. 
\end{enumerate}
\end{enumerate}
\end{theorem}

{\it Proof}. 1. Let $M\in\GWDO$.

(i)  $M$ {\em is a sum of modules} $V(I)$ {\em where} $I\in\CI(D_n,\gm)$, {\em see Proposition \ref{A15Apr17}.(1)}: The statement follows from
Proposition \ref{A15Apr17}.(1).

\medskip
(ii)  $\ann_{\mI_n}(M)=\ga_n$, by Lemma \ref{a15May19}  (since  $\ann_{\mI_n}(V(I))=\ga_n$). 

The statement (ii) means that $M$ is a $\Z^{n}$-graded $B_n$-module.
By Proposition \ref{A15Apr17}.(1), the functor $\GWO \rightarrow \LFm(D_n),\; M\mapsto M^{\gm}$ is an equivalence of categories with the inverse
$N\mapsto B_n\t_{D_n}N$.

2. Recall that $\gm=(H_1 -1, \ldots, H_m -1, H_{m+1},\ldots , H_l, H_{l+1}-\l_{l+1}, \ldots, H_n -\l_{n})$. Let $M\in\GWDO$.

(i)  $M$ {\em is a sum of modules} $V(I)$ {\em where} $I\in\CI(D_k, \gm')$, {\em see Proposition \ref{A15Apr17}.(2)}: The statement follows 
from Proposition \ref{A15Apr17}.(2). Since, for all $I\in\CI(D_k, \gm')$, 
$$V(I)^{\gm}=(P_m\t(B_k\t_{D_k}D_k/I))^{\gm}\simeq D_k/I,$$
\noindent the statement 2(a) follows from Proposition \ref{A15Apr17}.(2) and statement 1. Now, the statement 2(b) follows from the statement 2(a) and
statement 1.\;\;$\Box $

Let $A$ be an algebra and $M$ be an $A$-module. We denote by $[M]$ the isomorphism class of the $A$-module $M$ and $A-{\rm Mod}/\simeq\; $ 
is the set of all the isomorphism classes of $A$-modules. In particular, $\LFm (D_n)/\simeq\;$  is the set of   isomorphism classes of $D_n$-modules in $\LFm (D_n)$.
A category of modules is called a {\em category of finite representation type} if it contains only finitely many indecomposable modules up to isomorphism.
Definition of tame and wild category the reader can find in \cite{Drozd}. Notice that every category of finite representation type is tame but not vice versa.\\

{\bf Criterion for the category $\GWDO$ to be of finite representation type, tame or wild.} The next theorem is a criterion for the category $\GWDO$ to be of finite representation type, tame or wild.

\begin{theorem}\label{A14Apr17}
Let $\OO \in \cM_n/G$.
\begin{enumerate}
\item The category $\GWDO$ is of finite representation type iff $\OO = \Z^n$ and $\CD_\OO =\{ 1, \ldots , n\}$, 
and in this case the simple $\mI_n$-module $P_n=K[x_1,\ldots , x_n]$ is the unique indecomposable $\mI_n$-module in the category $\GWDO$. 
\item The category $\GWDO$ is tame  iff $|\CD_\OO |=n-1$, and in this case, 
up to order, $\OO =\Z^{n-1}\times (\l_n+\Z )$ for some $\l_n\in K$, $\CD_\OO=\{ 1, \ldots , n-1\}$ and $\{P_{n-1}\t M(i,\l)\, | \, i\in \N_+,$ where $\l= \l_n$ if $\l_n\not\in \Z$ and $\l =0$ if $\l_n\in \Z\}$ is the set of all indecomposable,
 pairwise non-isomorphic modules in $\GWDO$. 
\item The category $\GWDO$ is wild  iff $n\geq 2$ and $m:=|\CD_\OO |<n-1$, and 
in this case, up to order, $\CD_\OO=\{ 1, \ldots , m\}$ where $1\leq m<n-1$,  and $\{P_m\t (B_{n-m}\t_{D_{n-m}}N) \, | \, [N]\in \LF_{\gm'}(D_{n-m})/\simeq$ 
 and $N$ is an indecomposable $D_{n-m}$-module$\}$ 
 is a set of indecomposable, pairwise non-isomorphic modules in $\GWDO$ where $D_{n-m}=K[H_{m+1},\ldots , H_n]$ and $D_0=K$. 
\end{enumerate}
\end{theorem}

{\it Proof}. If $n\geq 2$ and $|\CD_{\O}|<n-1$ then by Theorem \ref{14Apr17} and \cite{Drozd-commutative} (see, also \cite{Ringel}), the category $\GWDO$ is wild.
If
either $n=1$ or $n\geq 2$ and $|\CD_{\O}|\geq n-1$ then by Theorem \ref{14Apr17} and \cite{Drozd-commutative} (see, also \cite{Ringel}), the category $\GWDO$ is tame.
Clearly, the category $\GWDO$ is of finite representation type iff $|\CD_{\O}|=n$.\;\;$\Box $

\begin{corollary}\label{b15Apr17}
 Let $\O\in\cM_n/G$. All modules in $\GWDO$ and $\WDO$ are
 equidimensional and the length of the module is equal the dimension of any 
 (generalized) weight component. In particular,
 all indecomposable generalized weight $\mI_n$-modules are 
 equidimensional.
\end{corollary}

{\it Proof}. The corollary follows from Theorem \ref{14Apr17}.\;\;$\Box $

The next corollary is a criterion for a generalized weight $\mI_n$-modules to be finitely generated.

\begin{corollary}\label{c15Apr17}
 Let $M$ be a generalized weight $\mI_n$-module.
 The $\mI_n$-module $M$ is finitely generated iff its support is a subset of a union of finitely
 many orbits in $\cM_n/G$ and the dimensions of all generalized weight components are bounded by a natural number.
\end{corollary}

{\it Proof}. The corollary follows from Theorem \ref{14Apr17} and Proposition \ref{A15Apr17}.\;\;$\Box $\\

{\bf Criterion for the category $\GWO$ to be of finite representation type, tame or wild.} Corollary \ref{a19May19} is
such a criterion.

\begin{corollary}\label{a19May19}
Let $\O\in\cM/G$.
\begin{enumerate}
 \item The category $\GWO$ is tame iff $n=1$.
 \item The category $\GWO$ is wild iff $n\geq 2$.
 \item None of the categories $\GWO$ is of finite representation type.
\end{enumerate}

\end{corollary}

{\it Proof}. The corollary follow from Theorem \ref{14Apr17}. \;\;$\Box $ \\

{\bf Explicit classes of indecomposable $\mI_n$-modules in $\GWO$.}
By Theorem \ref{14Apr17}, the problem of classifying indecomposable generalized weight $\mI_n$-modules in $\GWDO$ is
equivalent to the problem of classifying indecomposable modules in $\LFm(D_{n^{\prime}})$ for some $n^{\prime}\leq n$.
The set $\ind.\LFm(D_n)$ of isomorphism classes of indecomposable modules in $\LFm(D_n)$ is the union
$$\ind.\LFm(D_n)=\bigcup_{i\geq 1} \ind(D_n, \gm^{i})$$
\noindent where the set $\ind(D_n, \gm^{i})$ contains the isomorphism classes of all the indecomposable $D_n$-modules $M$
with $\gm^{i}M=0$. 
Clearly, 
$$\ind(D_n, \gm)=\{D_n/\gm\}\subseteq \ind(D_n, \gm^{2})\subseteq \ldots .$$
By Theorem of Drozd, see \cite{Drozd-commutative},

\medskip
\noindent  $\bullet\;\; \ind(D_n, \gm^{i})$  {\em is tame iff either $n=1$  or  $n=2$  and $m=1,2$}.

\medskip

{\bf Description of the set $\ind (D_2, \gm^{2})$.} Let $\L=D_2/\gm^{2}$,
$\gm=(h_1, h_2)$ where $h_1=H_1-\l_1$ and $h_2=H_2-\l_2$ for some $\l_1, \l_2\in K$ and $M\in \ind (D_2, \gm^{2})$.
Then $M=M_1\oplus M_2$ where $M_2=\gm M$ is a $\L$-module and $M_1$ is any (fixed) complement subspace of the vector
space $M_2$. Clearly, $M=0$ iff $M_2=M$ iff $M_1=0$.
The $\L$-module structure on $M$ is uniquely determined by the linear maps

$$\xymatrix{M_1 \ar@/^1pc/[r]^{h_1} \ar@/_1pc/[r]_{h_2} &M_2 },
\;\; m_1\mapsto h_1m_1, \;\;m_1\mapsto h_2m_1 \;\; ({\rm where}\;\; m_1\in M_1).$$
\smallskip

So, the problem of describing the set $\ind (D_2, \gm^{2})$ is `almost' equivalent to the problem of classifying 
indecomposable finite dimensional representations of the {\em Kronecker quiver}:

\[
\xymatrix{
1 \ar@/^1pc/[r]^{h_1} \ar@/_1pc/[r]_{h_2} &2 
}.\]

More precisely, every indecomposable finite dimensional representation of the Kronecker quiver $(M_1, M_2)$ such that
$M_1\neq 0$ belongs to $\ind (D_2, \gm^{2})$, and vice versa.
Up to isomorphism there are following $5$ series of indecomposable modules (in bracket bases of the vector spaces $M_1$ and $M_2$
are given):

\begin{enumerate}
 \item $K=D_2/\gm$.
 \item For each $n\geq 1,$ $M_1=\langle e_1, \ldots, e_n\rangle$, $M_1=\langle e_1^{\prime}, \ldots, e_{n+1}^{\prime}\rangle$,
 $h_1e_i=e_i^{\prime}$ and $h_2e_i=e_{i+1}^{\prime}$ for $i=1, \ldots, n$.
 \item For each $n\geq 1,$ $M_1=\langle e_1, \ldots, e_{n+1}\rangle$, $M_1=\langle e_1^{\prime}, \ldots, e_{n}^{\prime}\rangle$,
 $h_1e_i=e_i^{\prime}$ and $h_2e_{i+1}=e_{i}^{\prime}$ for $i=1, \ldots, n$, $h_1e_{n+1}=0$ and $h_2e_1=0$.
 \item For each $n\geq 1,$ $M_1=\langle e_1, \ldots, e_{n}\rangle$, $M_1=\langle e_1^{\prime}, \ldots, e_{n}^{\prime}\rangle$,
 $h_1e_i=e_i^{\prime}$ and $(h_2-\l)e_{i}=e_{i-1}^{\prime}$ for $i=1, \ldots, n$ where $e_0^{\prime}=0$ and $\l\in K$.
 \item For each $n\geq 1,$ $M_1=\langle e_1, \ldots, e_{n}\rangle$, $M_1=\langle e_1^{\prime}, \ldots, e_{n}^{\prime}\rangle$,
 $h_1e_i=e_{i-1}^{\prime}$ and $h_2e_{i}=e_{i}^{\prime}$ for $i=1, \ldots, n$ where $e_{0}^{\prime}=0$.
\end{enumerate}

The set $\ind.\LFm (D_n)$ is a disjoint union of subsets,

\begin{equation}\label{indUDn}
\ind.\LFm(D_n)=\bigsqcup_{I\in\CI(D_n,\gm)}\ind.\LFm(D_n, I)
\end{equation}

\noindent where the set $\ind.\LFm(D_n, I)$ contains all the indecomposable modules $M\in\ind.\LFm(D_n)$
with $\ann_{D_n}(M)=I$. By Theorem of Drozd, see \cite{Drozd-commutative},

\medskip
\noindent $\bullet\;\; \ind(D_n, I)$  {\em is tame iff either} $n=1$ {\em or} $n=2$ {\em and}
$I$ {\em contains a product  $h_1^{\prime}h_2^{\prime}$ of elements   $h_1^{\prime},h_2^{\prime}\in\gm$ such that their images in the
$K$-vector space $\gm/\gm^{2}$ are $K$-linearly independent (equivalently, are a basis)}.
\medskip

\noindent In the second case (i.e., $n=2$), the elements $h_1^{\prime}$ and $h_2^{\prime}$ are $K$-algebra generators for the algebra $\G$.
So, up to change of algebra generators, we can assume that $h_1h_2\in I$. Let $\G=D_2/(h_1h_2)=K[h_1, h_2]/(h_1h_2)$.
Then $\ind.\LFm(D_2, I)=\{M\in \ind_{f}(\G)\,|\,IM=0\}$ where $\ind_{f}(\G)$ is the set of isomorphism classes of indecomposable finite dimensional
left $\G$-modules.\\

{\bf Description of $\ind_{f}(\G)$}. 
Let $W=\langle h_1, h_2\rangle$ be a free (noncommutative) semigroup.
Each element (word) $w\in W$ is a unique product $w_1\cdots w_l$ where $w_i\in \{h_1, h_2\}$ and $l=1, 2, \ldots$.
The number $l=l(w)$ is called the {\em length } of the word $w$ and $W=\sqcup_{l\geq 1}W_l$, a disjoint union, where
$W_l$ is the set of all words of length $l$. 
The cyclic group of order $l$, $C_l=\langle \tau_l\rangle=\{\tau_l^{i}\,|\, i=0, \ldots, l-1\}$, where $\tau_l=(12\ldots l)$, 
acts on the set $W_l$ by the rule $\tau_l(w_1\cdots w_l)=w_{\tau_l(1)}\ldots w_{\tau_l(l)}$.
Let $W_l/C_l$ be the set of orbits. 
We say that two elements $w$ and $w^{\prime}$ of $W$ are {\em equivalent}, $w\sim w^{\prime}$, if they belong to the same orbit
($w\sim w^{\prime}$ iff $l(w)=l(w^{\prime})$ and $w=\tau_l^{i}(w^{\prime})$ for some $i$ where $l=l(w)$).
An orbit $\OOO\in W_l/C_l$ is called a {\em periodic} orbit if it contains an element $w$ such that $w=\theta^{i}$ for some
$\theta\in W$ and $i\geq 2$. 
We denote by $\CN$ the set of all {\em non-periodic orbits}.
The simple module $K=\G/(h_1, h_2)$  belongs to $\ind_f(\G)$.
The set of {\em non-simple} indecomposable finite dimensional $\G$-modules consists
of two sets of modules: the modules of the first and second type, see \cite{GP}:

\begin{equation}\label{ind12}
 \ind_f(\G)\setminus \{K\}=\ind_1(\G)\bigsqcup \ind_2(\G)
\end{equation}

\noindent where

\begin{enumerate}
 \item $ \ind_1(\G)=\{M_w\,|\,w\in W\}$ and $M_w=\langle e_1, e_2, \ldots, e_{l+1}\rangle$ where $l=l(w)$,
 $w=w_1\ldots w_l$ and $w_i\in\{h_1, h_2\}$,
 $$h_1e_i=\begin{cases}
e_{i+1}& \text{if } w_i=h_1,\\
0& \text{otherwise},\\
\end{cases}\;\;\;\;
h_2e_i=\begin{cases}
e_{i-1}& \text{if } w_{i-1}=h_2, i\geq 2,\\
0& \text{otherwise}.\\
\end{cases}$$
\item $\ind_2(\G)=\{N(\OOO, n, \l)\,|\, \OOO\in\CN, n\in\N_{+}, \l\in K^{*}\}$. 
Let $w=w_1\cdots w_l\in \OOO$ where $l=l(w)$. 
Then 
$$N(\OOO, n, \l)=\bigoplus_{i\in \Z/l\Z}N_i$$ 
\noindent is a direct sum of $n$-dimensional vector spaces
$N_i=K^{n}$ and the action of the elements $h_1$ and $h_2$ is given below.
Schematically, it can be represented by the following diagram

 \[
  \xymatrix{ N_1 \ar[r]^{\id} & N_2 \ar[r]^{\id} & \cdots  \ar[r]^{\id} &  N_l   \ar@/^2pc/[lll]^{J_n(\lambda)}  },
 \]

$$h_1|_{N_i}: N_i\to N_{i+1}, \;\; h_1|_{N_i}=\begin{cases}
\id & \text{if } i\neq l, w_i=h_1,\\
0& \text{if} i\neq l, w_i=h_2,\\
J_n(\l) & \text{if } i=l, w_l=h_1,\\
0 & \text{if } i=l, w_l=h_2,\\
\end{cases}$$

$$h_2|_{N_i}: N_i\to N_{i-1}, \;\; h_2|_{N_i}=\begin{cases}
0 & \text{if } i\neq l, w_i=h_1,\\
\id & \text{if} i\neq l, w_i=h_2,\\
0 & \text{if } i=l, w_l=h_1,\\
J_n(\l) & \text{if } i=l, w_l=h_2.\\
\end{cases}$$
 
\end{enumerate}

\noindent Up to isomorphism, the module $N(\OOO, n, \l)$ does not depend on the choice of the representative $w$ 
of the orbit $\OOO$.\\

{\bf Description of $\ind_f(A)$ where $A=K[h_1, h_2]/(h_1^{2}, h_2^{2})$}. The field $K$ is an algebraically closed field of
characteristic zero. Let $i=\sqrt{-1}$, $h_1^{\prime}=h_1+ih_2$ and $h_2^{\prime}=h_1-ih_2$.
Then $h_1^{\prime}h_2^{\prime}=h_1^{2}+h_2^{2}\in (h_1^{2}, h_2^{2})$, and so $A$ is tame, by Theorem of Drozd, see \cite{Drozd-commutative}.
Since the algebra $A$ is an epimorphic image of the algebra $\L=K[h_1, h_2]/(h_1^{2}, h_1h_2, h_2^{2})$, $\ind_f(A)=\ind_f(\L)\cup \{_A{}A\}$.

\begin{lemma}\label{a19Apr17}
 $\ind_f(A)=\ind_f(\L)\cup \{_A{}A\}.$
\end{lemma}

{\it Proof}. (i) $_A{}A$ is {\em indecomposable} (since $A$ is local).

\medskip
\noindent(ii) $_A{}A$ is {\em an injective module}: straightforward.
       
\medskip
\noindent(iii) {\em  Any finite dimensional $A$-module $M$ such that $\gm^{2}M\neq 0$ contain  $_A{}A$ where $\gm=(h_1, h_2)$}:
Since $\gm^{2}A\neq 0$, we can find a nonzero element $a\in A$ such that $\gm^{2}a\neq 0$.
Then $_A{}Aa\simeq A$, as required (since $\gm^{2}=(h_1h_2)$). 

\medskip
\noindent(iv) $\ind_f(A)=\ind_f(\L)\cup \{_A{}A\}$:
Let $M\in\ind_f(A)$. 
If $\gm^{2}M\neq 0$ then $M\simeq _A{}A$, by the statement (iii).
If $\gm^{2}M= 0$ then $M\in \ind_f(\L)$ since $\L=A/\gm^{2}$. \;\; $\Box$


\section{Generalized weight right $\mI_n$-modules}\label{GWRIM}


In this section, a classification of simple (generalized) weight  right $\mI_n$-modules
is given (Theorem \ref{RA11Apr17}). The category of weight right $\mI_n$-modules is a semisimple 
category (Theorem \ref{RB11Apr17}). An explicit description of generalized weight $\mI_n$-modules
is given (Theorem \ref{R14Apr17}).

The algebra $\mI_n$ admits an involution $\ast$ given by the rule, see \cite{Bav-gldim-intdif-JacAlg}: 
For $i=1, \ldots, n$, $$\der_i^{\ast}=\int_i, \;\; \int_i^{\ast}=\der_i\;\; {\rm  and}\;\; H_i^{\ast}=H_i.$$
Recall that an involution $\ast$ on $\mI_n$ is a $K$-algebra {\em anti-isomorphism} of $\mI_n$
$((ab)^{\ast}=b^{\ast}a^{\ast})$ such that $a^{\ast\ast}=a$ for all elements $a\in \mI_n^{\ast}.$
Clearly, the involution $\ast$ above acts as the identify map on the algebra $D_n$.

Every left $\mI_n$-module $M$ can be seen as a right $\mI_n$-module $M^{\ast}$ where
$M^{\ast}=M$, equality of vector spaces, and right $\mI_n$-module structure on $M$ is given by the rule:
For all $m\in M$ and $a\in \mI_n$, $ma:=a^{\ast}m$. Similarly, every right $\mI_n$-module $N$ can be seen
as a left $\mI_n$-module $N^{\ast}$ where $N^{\ast}=N$, equality of vector spaces, and, for all $n\in N$
and $a\in \mI_n$, $an:=na^{\ast}$. The functor 
$$\Md{\mI_n}\to {\rm Mod}\mbox{-}\mI_n,\;\; M\mapsto M^{\ast}$$
\noindent is an equivalence of categories with the inverse $N\mapsto N^{\ast}$.
Clearly, $M^{\ast\ast}=M$ and $N^{\ast\ast}=N$. 

\begin{example} Recall that the polynomial algebra $P_n$ is a left $\mI_n$-module isomorphic to the factor module
 $\mI_n/\mI_n(\der_1, \ldots, \der_n)\simeq K[\int_1, \ldots, \int_n]\overline{1}$ where
 $\overline{1}=1+\mI_n(\der_1, \ldots, \der_n)$ (since for all $\alpha\in \N^{n}$, $(\alpha ! )^{-1}\int^{\alpha}1=x^{\alpha})$. Hence, 
 \begin{equation}\label{Pn*}
  (P_n^{\ast})_{\mI_n}\simeq \mI_n/(\der_1, \ldots, \der_n)\mI_n\simeq\mI_n/(\int_1, \ldots, \int_n)\mI_n\simeq \tilde{1} K[\int_1^{\ast}, \ldots, \int_n^{\ast}]
  \simeq \tilde{1} K[\der_1, \ldots, \der_n]=\tilde{1} \CD_n\simeq \CD_n
 \end{equation}
 
 \noindent where $\tilde{1}=1+(\int_1, \ldots, \int_n)\mI_n$ and $\CD_n=K[\der_1, \ldots, \der_n]$ is a polynomial algebra.
 The algebra $\CD_n$ is a maximal commutative subalgebra of $\mI_n$.
 Let $e_1=(1,0, \ldots, 0)$, $\ldots, e_n=(0, \ldots, 0, 1)$ be the standard basis of the free abelian group $\Z^{n}=\oplus_{i=1}^{n}\Z e_i$.
 The right $\mI_n$-module $\CD_n=(P_n^{\ast})_{\mI_n}$ is simple (since $_{\mI_n}{P_n}$ is simple) and $\{\der^{\alpha}=\der_1^{\alpha_1}\cdots \der_n^{\alpha_n}\,|\,
 \alpha=(\alpha_1, \ldots, \alpha_n)\in\N^{n}\}$ is a $K$-basis of $\CD_n$.
 The right action of the generators $H_i, \der_i, \int_i$ ($i=1, \ldots, n$) of the algebra $\mI_n$ on $\der^{\alpha}$ are given below:
 $$\der^{\alpha} H_i=\der^{\alpha}(\alpha_i+1), \;\;\der^{\alpha}\der_i=\der^{\alpha +e_i}\;\; {\rm and}\;\; \der^{\alpha}\int_i=\begin{cases}
\der^{\alpha -e_i}& \text{if } \alpha_i\geq 1,\\
0& \text{if } \alpha_i=0.\\
\end{cases}$$

\end{example}

The definition of generalized weight right $\mI_n$-modules is given in the same way as their left counterparts. 
We add the subscript `r' to all the notation introduced for generalized weight left modules to indicate that we
deal with {\em right} modules.

Since the involution $\ast$ acts as the identity map on the polynomial algebra $D_n=K[H_1, \ldots, H_n]$,
we have, for each orbit $\O\in \cM_n/G$, 
\begin{equation}\label{RMMO}
\WO^{\ast}=\rWO,\;\;  \GWO^{\ast}=\rGWO, \;\; 
\end{equation}

\begin{equation}\label{RMMO1}
 \rWO^{\ast}=\WO,\;\;  \rGWO^{\ast}=\GWO,
\end{equation}

\begin{equation}\label{RMMO2}
 \rWIn=\bigoplus_{\O\in\cM_n/G}\rWO\;\;{\rm and}\;\; \rGWIn=\bigoplus_{\O\in\cM_n/G}\rGWO.
\end{equation}

So, for each $M\in\GWO$, $\Supp(M^{\ast})=\Supp(M)$ and $(M^{\ast})^{\gm}=M^{\gm}$ for all $\gm\in\Supp(M)$.

For each $\CD_{\O}, M(\CD_{\O})_{r}:=M(\CD_{\O})^{\ast}$ is a
simple right $\mI_n$-module with 
$$\Supp(M(\CD_{\O})_{r})= \Supp(M(\CD_{\O}))=S(\CD_{\O}),$$ 
\noindent see (\ref{SDCO}).
Then $\rGWDO:=\GWDO^{\ast}$ is the full subcategory of $\rGWO$
generated by the simple right weight $\mI_n$-module 
$M(\CD_{\O})_r$. There are precisely $2^{|\mD_{\O}|}$
such subcategories in the category $\rGWDO$, and 
\begin{equation}\label{rGWD}
\rGWO=\bigoplus_{\CD_{\O}\subseteq \mD_{\O}}\rGWDO,
\end{equation}
\noindent by Theorem \ref{15May19} (apply $^{\ast}$ to (\ref{GWDD})). \\

{\bf Description of simple weight modules.} We denote by $\widehat{\mI}_n({\rm weight})_r$ (resp., $\widehat{\mI}_n({\rm gen.\; weight})_r$)
the set of isomorphism classes of simple weight right (resp., generalized weight right) $\mI_n$-modules.
The next theorem classifies (up to isomorphism) all the simple weight right $\mI_n$-modules.

\begin{theorem}\label{RA11Apr17}
\begin{enumerate}
 \item $\widehat{\mI}_n(gen.\; weight)_r=\widehat{\mI}_n(weight)_r=\widehat{\mI}_1(weight)_r^{\otimes n}$, i.e.,
 any simple generalized weight right $\mI_n$-module is a simple weight right $\mI_n$-module, and vice versa;
 any simple weight right $\mI_n$-module $M$ is isomorphic to the tensor product $M_1\otimes \cdots \otimes M_n$
 of simple weight right $\mI_1$-modules and two such modules are isomorphic over $\mI_n$,
 $M_1\otimes \cdots \otimes M_n\simeq M_1^{\prime}\otimes \cdots \otimes M_n^{\prime},$
 iff for each $i=1, \ldots, n$, the $\mI_1$-modules $M_i$ and $M_i^{\prime}$ are isomorphic.
 Furthermore, $\widehat{\mI}_1(weight)_r=\{P_1^{\ast}\simeq \mI_1/\int\mI_1\simeq K[\der]$,
 $M(1,\l)^{\ast}\simeq B_1/(H-\l)B_1\,|\,\l\in K\}$, $\Supp(K[\der])=\N_{+}$ and $\Supp(M(1,\l)^{\ast})=\l+\Z.$

 \item For each simple weight right $\mI_n$-module $M= \bigotimes_{i=1}^{n}M_i$, $\Supp(M)=\prod_{i=1}^{n}\Supp(M_i).$
\end{enumerate}

\end{theorem}

{\it Proof}. The theorem follows at once from Theorem \ref{A11Apr17}, (\ref{RMMO}), (\ref{RMMO1})  and (\ref{RMMO2}).\;\;$\Box$

\begin{theorem}\label{RB11Apr17}
Every weight right $\mI_n$-module is a direct sum of simple weight right $\mI_n$-modules. In particular,
the category $\rWIn$ of weight right $\mI_n$-modules is a semisimple category.
\end{theorem}

{\it Proof}. The theorem follows from Theorem \ref{B11Apr17}, (\ref{RMMO}), (\ref{RMMO1}) and (\ref{RMMO2}).\;\;$\Box$\\

{\bf Explicit description of modules in $\rGWDO$.}

\begin{theorem}\label{R14Apr17}
Let $\OO \in \cM_n/G$.
\begin{enumerate}
\item Suppose that $\CD_\OO= \emptyset$,   $\gm\in \OO$   if $\mD_\O =\emptyset$ and $\gm =(H_1,\ldots , H_l, H_{l+1}-\l_{l+1}, \ldots, H_n-\l_n)\in S(\CD_\O )$ if $\mD_\OO = \{ 1, \ldots , l\}\neq \emptyset$, up to order. Then the functor
$$ \rGWDO \ra \LFm (D_n), \;\; M\mapsto M^\gm$$ is an equivalence of categories with the inverse 
$N\mapsto N\t_{D_n}B_n$, the induced functor. 
\item Suppose that $\CD_\OO\neq \emptyset$ and, up to order, $\CD_\OO = \{ 1, \ldots , m\}$, $\mD_\OO =\{ 1, \ldots , l\}$ for some $m$ such that $1\leq m \leq l\leq n$. Then $\mI_n=\mI_m\t \mI_{n-m}$ and 
$\OO = \Z^m\times \O'$ where $\O'=\Z^{l-m}\times\OO_1\times\cdots\times  \OO_{n-l}$ and $\OO_i\neq \Z$ for all $i=1,\ldots ,n-l$, and 
\begin{enumerate}
\item $\rGWDO = P_m^{\ast}\t \rGW (\mI_{n-m}, \CD_{\O'}'):=\{ P_m^{\ast}\t M\, | \, M\in \GW (\mI_{n-m},  \CD_{\O'}')\}$ with $ \CD_{\O'}'=\emptyset $.
\item Fix $\gm \in \OO$ such that $\gm =(H_1-1,\ldots , H_m-1, H_{m+1}, \ldots, H_l, H_{l+1}-\l_{l+1},\ldots , H_n-\l_n)$. 
Then $\gm' =(H_{m+1}, \ldots, H_l, H_{l+1}-\l_{l+1},\ldots , H_n-\l_n)\in\O'$ and  the functor 
$$ \rGWDO \ra \LF_{\gm'} (D_{n-m}), \;\; P_m^{\ast}\t M\mapsto M^{\gm'}$$
is an equivalence of categories with the inverse 
$N\mapsto P_m^{\ast}\t(N\t_{D_{n-m}}B_{n-m})$ where $D_{n-m}=K[H_{m+1},\ldots , H_n]$ and $D_0:=K$. 
\end{enumerate}
\end{enumerate}
\end{theorem}

{\it Proof}. The theorem follows from Theorem \ref{14Apr17} by
applying $^{\ast}$.\;\;$\Box$\\

{\bf Criterion for the category $\rGWDO$ to be of finite representation type, tame or wild.} Theorem \ref{RA14Apr17}
is such a criterion.

\begin{theorem}\label{RA14Apr17}
Let $\OO \in \cM_n/G$.
\begin{enumerate}
\item The category $\rGWDO$ is of finite representation type iff $\OO = \Z^n$ and $\CD_\OO =\{ 1, \ldots , n\}$, 
and in this case the simple right $\mI_n$-module $P_n^{\ast}$ is the unique indecomposable $\mI_n$-module in the category $\rGWDO$. 
\item The category $\rGWDO$ is tame  iff $|\CD_\OO |=n-1$, and in this case, 
up to order, $\OO =\Z^{n-1}\times (\l_n+\Z )$ for some $\l_n\in K$, $\CD_\OO=\{ 1, \ldots , n-1\}$ and $\{P_{n-1}^{\ast}\t M(i,\l)\, | \, i\in \N_+,$ where $\l= \l_n$ if $\l_n\not\in \Z$ and $\l =0$ if $\l_n\in \Z\}$ is the set of all indecomposable,
 pairwise non-isomorphic modules in $\rGWDO$. 
\item The category $\rGWDO$ is wild  iff $n\geq 2$ and $m:=|\CD_\OO |<n-1$, and 
in this case, up to order, $\CD_\OO=\{ 1, \ldots , m\}$ where $1\leq m<n-1$,  and $\{P_m^{\ast}\t (N\t_{D_{n-m}}B_{n-m}) \, | \, [N]\in \LF_{\gm'}(D_{n-m})/\simeq$ 
 and $N$ is an indecomposable $D_{n-m}$-module$\}$ 
 is a set of indecomposable, pairwise non-isomorphic modules in $\rGWDO$ where $D_{n-m}=K[H_{m+1},\ldots , H_n]$ and $D_0=K$. 
\end{enumerate}
\end{theorem}

{\it Proof}. The theorem follows from Theorem \ref{A14Apr17}  by
applying $^{\ast}$.\;\;$\Box$

\begin{corollary}\label{Rb15Apr17}
 Let $\O\in\cM_n/G$. All modules in $\rGWDO$ and $\rWDO$ are
 equidimensional and the length of the module is equal the dimension of any of
 (generalized) weight component. In particular,
 all indecomposable right generalized weight $\mI_n$-modules are 
 equidimensional.
\end{corollary}

{\it Proof}. The corollary follows from Corollary \ref{b15Apr17}.\;\;$\Box $

The next corollary is a criterion for a generalized weight right $\mI_n$-modules to be finitely generated.

\begin{corollary}\label{Rc15Apr17}
 Let $M$ be a generalized weight right $\mI_n$-module.
 The $\mI_n$-module $M$ is finitely generated iff its support is a subset of a union of finitely
 many orbits in $\cM_n/G$ and the dimension of all generalized weight components are restricted by a natural number.
\end{corollary}

{\it Proof}. The corollary follows from Corollary \ref{c15Apr17}.\;\;$\Box $ 

\begin{corollary}\label{ra21May19}
 \begin{enumerate}
  \item Every module $M\in\rGWO$ is a unique direct sum of
  absolutely prime generalized weight right $\mI_n$-modules, 
  and this direct sum is $M=\bigoplus_{\CD_{\O}\subseteq \mD_{\O}}M_{\CD_{\O}}$ where $M_{\CD_{\O}}\in \rGWDO.$
  \item Every generalized weight right $\mI_n$-module is a unique sum of
  absolutely prime generalized weight right $\mI_n$-modules.
 \end{enumerate}

 \end{corollary}

 {\it Proof}. The corollary follows 
 from Corollary \ref{a21May19} by applying $^{\ast}$. \;\;$\Box$\\

{\bf Criterion for the category $\rGWO$ to be of finite representation type, tame or wild.} Corollary \ref{ra19May19} is
such a criterion.

\begin{corollary}\label{ra19May19}
Let $\O\in\cM/G$.
\begin{enumerate}
 \item The category $\rGWO$ is tame iff $n=1$.
 \item The category $\rGWO$ is wild iff $n\geq 2$.
 \item None of the categories $\rGWO$ is of finite representation type.
\end{enumerate}

\end{corollary}

{\it Proof}. The corollary follow from Corollary \ref{a19May19}
by
applying $^{\ast}$. \;\;$\Box $ \\

Using the involution, we can consider right analogues of indecomposable $\mI_n$-modules considered at the end of
Section \ref{CSGWM}. 
We leave this to the interested reader.

$${\rm\bf Acknowledgment.}$$

The first author  is partly supported by Fapesp grant (2017/02946-0) and by a grant of the University of 
Sheffield.
The third  author is partly supported by CNPq grant (301320/2013-6) and by Fapesp grant (2014/09310-5).
This work was done during the visit of the first author to the University of S\~{a}o Paulo and the University of 
Belo Horizonte (UFMG) whose hospitality and support are greatly acknowledged.

\small{







}


\begin{thebibliography}{99}





\bibitem{Bav-FinDimExt-91} V. V. Bavula, 
Finite-dimensionality of ${\rm Ext}_n$ and ${\rm Tor}_n$ of simple modules over a class of algebras, 
{\em Funct. Anal. Appl.}  {\bf  25} (1991) no. 3, 229--230.

\bibitem{Bav-GWA-1992} V. V. Bavula, 
Generalized Weyl algebras and their representations, 
{\em Algebra i Analiz}  {\bf 4} (1992), no. 1, 75--97; 
English transl. in {\em St. Petersburg Math. J.} {\bf 4} (1993) no. 1, 71--92.

\bibitem{Bav-92} V. V. Bavula, 
Simple $D[X,Y;\sigma, a]$-modules,
{\em Ukra\"{\i}n. Mat. Zh.} {\bf 44} (1992), no. 12, 1628--1644;
English transl. in {\em Ukrainian Math. J.} {\bf 44} (1992),
no. 12, 1500--1511.


\bibitem{Bav-Schurian} V. V. Bavula, 
Each Schurian algebra is tensor-simple, 
{\em Comm. Algebra} {\bf  23} (1995) no. 4, 1363--1367.


\bibitem{Bav-Bek}  V. V. Bavula and V. Bekkert,
 Indecomposable representations of generalized Weyl algebras, 
 {\em Comm. Algebra} {\bf  28} (2000) no. 11, 5067--5100.
 
\bibitem{BBF-1} V. V. Bavula, V.  Bekkert  and V. Futorny,
Indecomposable generalized weight  modules over the algebra of polynomial
integro-differential operators,
{\em Proc. Am. Math. Soc.} {\bf 146} (2018) no. 6, 2373--2380.

\bibitem{algintdif}  V. V. Bavula,    
The algebra of integro-differential operators on a polynomial algebra, 
{\em J.  Lond. Math. Soc.  (2)} {\bf 83} (2011) no. 2, 517--543.



\bibitem{intdifaut}    V. V. Bavula,  
The group of automorphisms of the algebra of polynomial integro-differential operators,
{\em  J. Algebra}  {\bf  348} (2011) 233--263.

\bibitem{algintdifline}  V. V. Bavula,    
The algebra of integro-differential operators on an affine line and its modules, 
{\em J. Pure Appl. Algebra} {\bf  217} (2013) no. 3, 495--529.

 \bibitem{indtif-bimod}  V. V. Bavula, 
 The algebra of polynomial integro-differential operators is a holonomic bimodule over the subalgebra of polynomial differential operators,
 {\em Algebr. Represent. Theory} {\bf  17} (2014) no. 1, 275--288.
 
 \bibitem{Bav-Lu-2016-Isr} V. V. Bavula and T. Lu, 
 The quantum Euclidean algebra and its prime spectrum,  
 {\em Israel J. Math.}, {\bf 219} (2017) no. 2, 929--958. 

\bibitem{Bav-gldim-intdif-JacAlg} V. V. Bavula, 
The global dimension of the  algebras of polynomial integro-differential operators $\mI_n$ and the Jacobian algebras $\mA_n$, {\em J. Algebra Appl.}, 
(2020) 2050030 (28 pages), 
DOI: 10.1142/S0219498820500309, 
arXiv: 1705.05227.



\bibitem{Bav-Lu-N3} V. V. Bavula and  T. Lu, Torsion simple modules over the quantum spatial ageing algebra, {\em Comm. Algebra}, {\bf 45} (2017) no. 10, 4166-4189.


\bibitem{Bav-Lu-N8} V. V. Bavula and  T. Lu, Prime ideals
of the enveloping algebra of the Euclidean algebra and a classification
of its simple weight modules, {\em J. Math. Phys.}, {\bf 58} (2017) no. 1,  011701, 33p,  DOI: 10.1063/1.4973378.


\bibitem{Bav-Lu-N9} V. V. Bavula and  T. Lu, Classification of simple weight modules over  the Schr\"{o}dinger algebra, {\em Canad. Math.  Bull.}, {\bf 61} (2018), no. 1, 16-39. 



\bibitem{Bav-Lu-N5} V. V. Bavula and  T. Lu,  The universal enveloping algebra $U(\sl_2\ltimes V_2)$, its prime spectrum and a classification of its simple weight modules, {\it J. Lie Theory},    {\bf 28}  (2018)  no. 2, 525--560.  



\bibitem{Bav-Lu-N7} V. V. Bavula and  T. Lu, The prime spectrum of the algebra $\mathbb{K}_q[X,Y]\rtimes U_q(\sl_2)$ and a classification of simple weight modules, {\it J. Noncomm. Geometry}, {\bf 12} (2018), no. 3, 889--946. 

\bibitem{Bek-Ben-Fut-2004} V. Bekkert, G. Benkart and V. Futorny, 
Weight modules for Weyl algebras. {\em  Kac-Moody Lie algebras and related topics}, 17--42, Contemp. Math., 343, Amer. Math. Soc., Providence, RI, 2004.


\bibitem{Block-IrrRepsl2} R. E. Block, 
The irreducible representations of the Lie algebra $\mathfrak{sl}(2)$ and of the Weyl algebra, 
{\em Adv. Math.} {\bf 39} (1981) 69--110.

\bibitem{Drozd-commutative} Yu. A. Drozd,
Representations of commutative algebras, 
{\em Funct. Anal. Appl.} {\bf 6} (1972) 286--288.

\bibitem{Drozd} Yu. A. Drozd, 
Tame and wild matrix problems,
{\em Representation theory, II (Proc. Second Internat. Conf., Carleton Univ., Ottawa, Ont., 1979)}, pp. 242--258, Lecture Notes in Math., 832, Springer, Berlin-New York, 1980.


\bibitem{Fut-Gr-Maz-2014}  V. Futorny, D. Grantcharov and  V. Mazorchuk, 
Weight modules over infinite dimensional Weyl algebras,
{\em Proc. Am. Math. Soc.} {\bf 142} (2014) no. 9, 3049--3057.

\bibitem{Fut-Iyer-2016} V.  Futorny and U. Iyer,  
Representations of $D_q(k[x])$,
{\em Israel J. Math.} {\bf 212} (2016) no. 1, 473--506.

\bibitem{GP} I. M. Gelfand and V. A. Ponomarev, 
Indecomposable representations of the Lorenz group, 
{\em Russian Math. Surveys} {\bf 23 } (1968) 1--58.

 \bibitem{Guo-Regen-Ros-2014} L. Guo,  G. Regensburger and  M. Rosenkranz, 
On integro-differential algebra, 
{\em J. Pure Appl. Algebra} {\bf 218} (2014) 456--471.


\bibitem{Hartwig-2006}  J. Hartwig, 
Locally finite simple weight modules over twisted generalized Weyl algebras, 
{\em J. Algebra} {\bf 303} (2006) no. 1, 42--76.

\bibitem{Hartwig-2011} J. Hartwig, 
Pseudo-unitarizable weight modules over generalized Weyl algebras, 
{\em J. Pure Appl. Algebra} {\bf  215}  (2011) no. 10, 2352--2377.


\bibitem{Maz-Tur-1999} V. Mazorchuk and L. Turowska, 
Simple weight modules over twisted generalized Weyl algebras, 
{\em Commun. Algebra} {\bf 27} (1999) no. 6, 2613--2625.

\bibitem{Lu-Maz-Zha0-2015} R.  Lu, V.  Mazorchuk and K. Zhao, 
Simple weight modules over weak generalized Weyl algebras, 
{\em  J. Pure Appl. Algebra} {\bf 219}  (2015) no. 8, 3427--3444.


\bibitem{MR}  J. C.  McConnell and  J. C. Robson,
 Homomorphisms and extensions of modules over certain differential polynomial rings, 
 {\em J. Algebra} {\bf  26} (1973) 319--342.
 
 \bibitem{Ringel}  C.~M.~Ringel,
The representation type of local algebras.
{\em Proceedings of the International Conference on Representations of Algebras (Carleton Univ., Ottawa, Ont., 1974),
 Paper No. 22}, 24 pp. Carleton Math. Lecture Notes, No. 9, Carleton Univ., Ottawa, Ont., 1974.

\bibitem{Ship-2010} I. Shipman, 
Generalized Weyl algebras: category O and graded Morita equivalence, 
{\em J. Algebra} {\bf 323} (2010) no. 9, 2449--2468.











\end{thebibliography}
\end{document}